\documentclass{article}
\usepackage[colorlinks=true]{hyperref}
\usepackage[T1]{fontenc}
\usepackage[english]{babel}
\usepackage{amsmath}
\usepackage{amssymb}
\usepackage{amsfonts}
\usepackage{latexsym}
\usepackage{epsfig}
\usepackage{color}
\usepackage{color}
\normalbaselines
\newtheorem{theorem}{Theorem}[section]

\newtheorem{proposition}[theorem]{Proposition}

\newtheorem{remark}[theorem]{Remark}
\newtheorem{example}[theorem]{Example}

\newenvironment{Proof}{\removelastskip\par
\noindent{\em Proof.} \rm}{\penalty-20\null\hfill$\square$\par\medbreak}

\newcommand{\w}{{{\omega} }}

\title{\bf Exact boundary controllability  and stabilizability  of a degenerated Timoshenko beam \footnote{This work has been supported by the DFG Funding (grants No. 504042427), the National Natural Science Foundation of China (grants No. 62588101 and No. 12271035) and Beijing Municipal Natural Science Foundation (grant No. 1232018)}}

\author{
G\"unter Leugering\thanks{Department Mathematik, LS Angewandte Mathematik 2,
Friedrich-Alexander-Universit\"at Erlangen-N\"urnberg, 
D-91058 Erlangen 
(Germany)} \and Yue Wang * \thanks{* Corresponding Authour (yue.wang@fudan.edu.cn) Center for Applied Mathematics, Fudan University, 200433 Shanghai (China)}
\and Qiong Zhang\thanks{ 
School of Mathematics and Statistics, Beijing Institute of Technology, 100081 Beijing (China)
}
}

\begin{document}
\maketitle
\begin{abstract}
This paper investigates the boundary controllability and stabilizability of a Timoshenko beam subject to degeneracy at one end, while control is applied at the opposite boundary. Degeneracy in this context is measured by the real parameters for $\mu_a\in [0,2)$ for $a\in\{K,EI\}$, where $K(x)$ denotes shear stiffness and $EI(x)$ bending stiffness. We differentiate between weak degeneracy $\mu_a\in [0,1)$ and strong degeneracy $\mu_a\in [1,2)$, which may occur independently in shear and bending. Our study establishes observability inequalities for both weakly and strongly degenerate equations under Dirichlet, Robin, and Neumann boundary conditions. Using energy multiplier techniques and the Hilbert Uniqueness Method (HUM), we derive conditions for exact boundary controllability and show that appropriate boundary state and velocity feedback controls at the non-degenerate end can stabilize the system exponentially. Extending results previously obtained for the 1-dimensional wave equation in \cite{AlabauCannarsaLeugering2017}, this study highlights new control strategies and stabilization effects specific to the degenerate Timoshenko beam system, addressing challenges pertinent to real-world structural damping and control applications.

\end{abstract}

\section{Introduction}
The boundary controllability and stabilization of beam structures are essential considerations in fields ranging from civil engineering to aerospace applications, where precise control over vibrations, deflections, and stability are critical for performance and safety. The Timoshenko beam model, which accounts for both shear deformation and rotational inertia, provides a robust framework for analyzing these factors. This model is particularly relevant for slender beams or structures requiring high precision, where simplified assumptions may not capture the nuanced behaviors due to shear and rotational effects. In recent years, the analysis of beams with varying material properties has gained attention, particularly in cases where structural degradation or material inhomogeneities lead to degeneracies at the boundaries. Degeneracy, as modeled here, impacts energy dynamics and presents unique challenges for control and stabilization, as standard methods often rely on uniform stiffness and density assumptions.

In this study, we investigate the Timoshenko beam with degeneracy at one boundary, characterized by vanishing stiffness parameters $K(x)$ or $EI(x)$ near the degenerate end. This problem extends classical studies on beam controllability and stabilization, such as those in \cite{AlabauCannarsaLeugering2017, Fragnelli_etal_2024, LagneseLeugeringSchmidt1993b}, by introducing boundary conditions that account for both weak and strong degeneracy. We consider a straight linear planar Timoshenko beam, where $w(t,x)$ represents the vertical deflection and $\psi(t,x)$ the shear angle, with $(t,x)\in (0,T)\times (0,\ell)$. The beam's physical parameters are $K(x)=G(x)h$ and $EI(x)$, where $G(x)$ the shear modulus, $h$ the thickness, $E(x)$ Young's modulus, and $I$ the moment of the cross section around the axis sticking out of the plane. The line density $\rho=m_0h$ and rotational density $I_\rho=m_0I$ as constants with given mass $m_0$. In this article, we keep the thickness constant and possibly small. Following \cite{AlabauCannarsaLeugering2017}, we apply boundary conditions at the degenerate end $x=0$ based on weak or strong degeneracy, parametrized by $\mu_K, \mu_{EI}, K,EI\in \mathcal{A}$ as defined in \eqref{eq:hp_a}, with controls applied at $x=\ell$. The governing equations for the system are
\begin{align}\label{Timo_weak_strong_Neumann}
& \rho w_{tt}-\big(K(x)(w_x+\psi)\big)_x=0,\; (t,x)\in Q,\notag\\
& I_\rho \psi_{tt}-(EI(x)\psi_x)_x+K(x)(w_x+\psi)=0,\; (t,x)\in Q,\notag\\
& w(t,0)=0, \; \text{ if } \mu_K\in[0,1),\; \lim\limits_{x\rightarrow 0^+}K(x)(w_x+\psi)(t,x)=0\;  \text{ if } \; \mu_K\in [1,2),\; t\in (0,T),\\
& \psi(t,0)=0, \; \text{ if } \mu_{EI}\in[0,1),\; \lim\limits_{x\rightarrow 0^+}EI(x)\psi_x(t,x)=0\;  \text{ if } \; \mu_{EI}\in [1,2),\; t\in (0,T),\notag\\
& K(\ell)(w_x(t,\ell)+\psi(t,\ell))+\gamma w(t,\ell)=f_1(t), \; EI(\ell)\psi_x(t,\ell)+\delta \psi(t,\ell)=f_2(t), \; t\in (0,T),\notag\\
& w(0,x)=w^0(x), \; w_t(0,x)=w^1(x), \; \psi(0,x)=\psi^0(x), \; \psi_t(0,x)=\psi^1(x), \; x\in (0,\ell),\notag
\end{align}
where $Q:=(0,T)\times(0,\ell),\; \gamma,\;\delta\geq 0$. At the end $x=\ell$, for $\gamma,\delta>0$, we have Robin controls; while for $\gamma=\delta=0$, \eqref{Timo_weak_strong_Neumann}$_5$ reflects Neumann-controls. Dirichlet controls are also considered by replacing \eqref{Timo_weak_strong_Neumann}$_5$ with $w(t,\ell)=f_1(t), \psi(t,\ell)=f_2(t)$. In the sequel, we refer to \eqref{Timo_weak_strong_Neumann} as the standard model, where we change eventually \eqref{Timo_weak_strong_Neumann}$_{5,6}$ if different boundary and initial conditions apply. This problem extends classical studies on beam controllability and stabilization, such as those in \cite{LagneseLeugeringSchmidt1993b},  by introducing boundary conditions that account for both weak and strong degeneracy. These conditions, which are determined by the parameters $\mu_K, \mu_{EI}$, enable us to examine scenarios where degeneracy in bending or shear may occur independently. The study thereby closes the gap between degenerate string systems and more complex beam structures and contributes to the ongoing development of control methods for non-uniform and damaged structures.

In the non-degenerate case, one has $K(x)\geq k_0>0, EI(x)\geq \gamma_0>0$. In the case of static degeneration or static damage, the case that we discuss in this article, the shear modulus  $G(x)$ and Young's modulus  $E(x)$ may vanish at $x=0$, while all other physical parameters remain constant. As longitudinal displacements decouple from $w,\psi$ for a straight linear beam, we omit the corresponding equation that would have to be added to \eqref {Timo_weak_strong_Neumann}. The degrees of degeneracy for $K$ and $EI$ may differ, which makes parameter relationships critical, so we retain the original notation without dimensionless reduction. Also, to avoid vanishing $\rho$, thickness is held constant, with $K(x)=G(x)h, EI(x):=E(x)I$; for the sake of simplicity, we write $a\in \{K,EI\}$ to express that $a(x)$ is either $K(x)$ or $EI(x)$.

Weak degeneracy occurs for $\mu_a \in [0,1)$ and strong degeneracy for $\mu_a\in [1,2)$ as defined in \eqref{eq:hp_a}, applying separately to $K$ and $EI$.
We also refer to a very recent result in \cite{Fragnelli_etal_2024}, where the degenerated Euler-Bernoulli beam is treated with respect to controllability and stabilizability. Therefore, the current article closes the gap between degenerated strings and Euler-Bernoulli beams. It should be noted that, in principle, one can retrieve the Euler-Bernoulli beam equation from the Timoshenko beam equation if one lets the shear stiffness tend to infinity. The analysis in this article, however, does not carry over to the Euler-Bernoulli beam.

\bigskip

As already stated, we may also consider Dirichlet controls at $x=\ell$, modifying system \eqref{Timo_weak_strong_Neumann} by replacing \eqref{Timo_weak_strong_Neumann}$_5$ with
\begin{align}\label{Timo_weak_strong_Dirichlet}
 w(t,\ell)=f_1(t), \; \psi(t,\ell)=f_2(t), \; t\in (0,T).
\end{align}

Based on the standard model \eqref{Timo_weak_strong_Neumann} with Neumann/Robin controls, or \eqref{Timo_weak_strong_Neumann} with Dirichlet controls \eqref{Timo_weak_strong_Dirichlet}, we now formulate our controllability and observability problems. For the problem of \textbf{exact boundary controllability}, we seek control functions $f_1,f_2\in U$ such that the corresponding solution $(w,\psi)\in C(0,T;X\times Y)$ of \eqref{Timo_weak_strong_Neumann} (or \eqref{Timo_weak_strong_Neumann} with Dirichlet controls \eqref{Timo_weak_strong_Dirichlet}) satisfies
\begin{align}\label{end}
w(T,x)=w_T^0(x),\; w_t(T,x)=w_T^1(x); \psi(T,x)=\psi^0_T(x), \; \psi_t(T,x)=\psi^1_T(x), \; x\in (0,\ell),
\end{align}
where $w^0,\psi^0, w^0_T, \psi^0_T\in X$ and $ w^1,\psi^1, w^1_T, \psi^1_T\in Y$, with spaces $U,X,Y$ suitably chosen. Due to the system's time reversibility, we consider null controllability, setting $w^0_T=w^1_T=\psi_T^0=\psi^1_T=0$ without loss of generality, leading to the state constraints studied in this article as
\begin{align}\label{null_end}
w(T,x)=0,\; w_t(T,x)=0,\; \psi(T,x)=0, \; \psi_t(T,x)=0, \; x\in (0,\ell).
\end{align}
As for the problem of \textbf{exact boundary observability}, we examine the homogeneous system where the boundary condition satisfies
\begin{align}\label{Timo_weak_strong_Neumann_homogen}
& K(\ell)(w_x(t,\ell)+\psi(t,\ell))+\gamma w(t,\ell)=0, \; EI(\ell)\psi_x(t,\ell)+\delta \psi(t,\ell)=0, \; t\in (0,T),
\end{align}
with $\gamma,\delta\geq 0$, corresponding to the Robin/Neumann-control problem \eqref{Timo_weak_strong_Neumann} or for Dirichlet conditions instead of \eqref{Timo_weak_strong_Neumann}$_5$
\begin{align}\label{Timo_weak_strong_Dirichlet_homogen}
& w(t,\ell)=0, \; \psi(t,\ell)=0, \; t\in (0,T).
\end{align}
 Observability then means that out of given extra data on the boundary - such as $w(t,\ell), \psi(t,\ell)$ in Neumann case or $K(\ell)(w_x+\psi)(t,\ell), EI(\ell)\psi_x(t,\ell)$ in the Dirichlet case - we are able to recover the initial data and, consequently, the full solution. We will state the problems of exact boundary controllability and observability more precisely, when we have clarified the space set-up in the coming sections. We note that non-degenerate cases have been studied extensively, e.g. in \cite{KimRenardy1987} for the stabilization of the constant coefficient model, and \cite{Ammar-Khodja_etal2007} for varying coefficients with essentially one velocity feedback in the shear angle at both ends. In \cite{Alabau2007}, the author considers a constant coefficient Timoshenko beam, but with only one, in fact, non-linear distributed feedback. Problems of stabilization of Timoshenko beams in various contexts have been investigated until very recently. See e.g. \cite{Krstic_etal2024} for the state of the art. An in-depth analysis of dynamic non-uniform Timoshenko beams has been provided in \cite{Belishev2010} including the construction of fundamental solutions. Controllability properties have been investigated by \cite{Taylor00}, by the method of extending the initial data to a semi-infinite interval and the use of locals smoothing, by
\cite{LagneseLeugeringSchmidt1993a}, for planar networks of Timoshenko beams with varying coefficients by the HUM Method, \cite{LagneseLeugeringSchmidt1993b} for such networks using characteristics.

The article is organized as follows: Section \ref{preliminaries}  introduces properties of the degenerate coefficients and suitable weighted Sobolev spaces for formulating the problem mathematically. Section \ref{multiplier_identities}  provides essential multiplier identities, followed by Section \ref{inequalities}, where we establish direct and so-called inverse inequalities crucial for well-posedness under homogeneous boundary conditions \eqref{Timo_weak_strong_Neumann_homogen} (or \eqref{Timo_weak_strong_Dirichlet_homogen} for Dirichlet-case), which, in turn, is treated in Section \ref{well-posendess} along with controlled systems. Sections  \ref{observability} and \ref{controllability} present and prove observability and controllability results for Dirichlet, Robin, and Neumann controls.  Lastly, Section \ref{stabilizability} addresses boundary exponential stabilizability of the Timoshenko system.  In all sections, we will discriminate between weak and strong degeneration. This work extends existing theory by rigorously addressing exact controllability and exponential stabilizability of Timoshenko beams with boundary degeneracies, offering new insights into the control of degenerate elastic systems.


\section{Assumptions and preliminaries}\label{preliminaries}
\subsection{Assumptions and properties of the degenerate coefficients}
We say a function $a\in\mathcal C([0,\ell])\cap \mathcal C^1(]0,\ell ])$ is in the class $\mathcal{A}$ if it  satisfies:
\begin{equation}\label{eq:hp_a}
\begin{cases}
 {\rm (i)} &  a(x)>0\;\;\forall x\in ]0,\ell ]\,,\;\;
a(0)=0\,,
\\
 {\rm (ii)}& \mu_a:=\sup_{0<x\leqslant \ell}\dfrac{x|a'(x)|}{a(x)}<2 \,,
\\
 {\rm (iii)}&
a\in \mathcal C^{[\mu_a]}([0,\ell ]),
\end{cases}
\end{equation}
where $[\cdot]$ stands for the integer part. In the sequel, we make the assumption
\begin{equation}\label{H}
 {\bf (H)}:   \quad   \quad    K,\; EI \in\mathcal{A}.
\end{equation}

\begin{proposition}\label{re:a}
We recall from \cite{AlabauCannarsaLeugering2017}  some useful properties of $a$ in the class $\mathcal A$, adapted to the interval $[0,\ell]$.
\begin{enumerate}
 \item[{\rm (i)}] We have
\begin{equation}\label{eq:a}
a(x)\ell^{\mu_a}\geqslant a(\ell ) x^{\mu_a},\qquad\forall x\in [0,\ell ]\,.
\end{equation}
Hence,  $1/a\in L^1(0,\ell )$ if $\mu_a\in [0,1)$.
\item[{\rm (ii)}] For  all $\mu_a\in[0,2)$,  the inequality $a(\ell)x^2\leq \ell^2 a(x)$ holds; and for $\mu_a\in [0,1]$, we have in particular $a(\ell)x\leq \ell a(x)$.
\item[{\rm (iii)}] When $\mu_a\in [1,2)$, it follows that $1/a\notin L^1(0,\ell )$.
\end{enumerate}
\end{proposition}

Following the approach in \cite{AlabauCannarsaLeugering2017}, for $a\in {\cal A}$, we  can define the weighted Sobolev space
$V^1_{a}(0,\ell) := \{u \in L^2(0,\ell )\;|\;  u$ is  locally absolutely continuous in $ ]0,\ell ], $
 and $ au' \in L^2(0,\ell )\}$,  equipped with the norm  $(\|\sqrt{a}u\|^2_{L^2(0,\ell )}+u^2(\ell))^{1\over2}$ for
 $u\in V^1_{a}(0,\ell)$.   In the special case where   $u(\ell)=0$, the equivalent norm is given by  $\|\sqrt{a}u\|_{L^2(0,\ell )}$.

\subsection{Weighted Sobolev Spaces for Degenerate Systems}
Following the approach in \cite{AlabauCannarsaLeugering2017}, we introduce the fundamental weighted Sobolev space associated with degenerate system \eqref {Timo_weak_strong_Neumann}. Let $V^1_{K,EI}(0,\ell )$ denote the space of function pairs $(w,\psi)\in L^2(0,\ell )^2$ satisfying
\begin{equation}\label{eq:H1a}
\begin{cases}
{\rm (i)}&
w, \psi\;\;\mbox{are  locally absolutely continuous in}\;]0,\ell ],
\;\mbox{and}
\\
{\rm (ii)}&
\sqrt{K}(w'+\psi), \sqrt{EI}\psi' \in L^2(0,\ell ).
\end{cases}
\end{equation}
It is straightforward to verify that $V^1_{K,EI}(0,\ell )$ is a Hilbert space with the scalar product
\begin{align}\label{scalar_product}
\langle (w,\psi),(z,\xi)\rangle_{1,K,EI}:=&\int_0^\ell\big[K(x)(w'(x)+\psi(x))(z'(x)+\xi(x))
\\& +EI(x)\psi'(x)\xi'(x)+ w(x)z(x)+\psi(x)\xi(x)\big]dx
\end{align}
for all $(w,\psi),(z,\xi)\in V^1_{K,EI}(0,\ell )$
and the associated norm.
Let us also introduce the bilinear form
\begin{align}\label{bilinear_form}
& B\;: \;V_{K,EI}^1(0,\ell)\times V_{K,EI}^1(0,\ell)\;  \rightarrow \; \mathbb{R},\\
&B((w,\psi),(z,\xi)):= \int\limits_0^\ell\big[ K(x)(w'(x)+\psi(x))(z'(x)+\xi(x))+EI(x)\psi'(x)\xi'(x)\big] dx.
\end{align}
The corresponding semi-norm is defined by $$|(w,\psi)|_{1,K,EI}:=\sqrt{B\big((w,\psi),(w,\psi)\big)}.$$
 This semi-norm $|\cdot|_{1,K,EI}$ is an equivalent norm on the closed subspace of $V^1_{K,EI,0}(0,\ell )$ defined as $$V^1_{K,EI,0} (0,\ell )=\big\{(w,\psi) \in V^1_{K,EI}(0,\ell ) \;\big| \; w(\ell )=\psi(\ell)= 0\big\}.$$
  This fact is a simple consequence of the following version of Poincar\'e's inequality.
\begin{proposition}\label{Poincare_inequality}
 Assume \eqref{eq:hp_a} and $\mu:=\max\{\mu_K,\mu_{EI}\}$  for simplicity.
 Then
\begin{equation}\label{eq:poincare}
\|(w,\psi)\|^2_{L^2(0,\ell)^2}\leq C_{D,K,EI}\; |(w,\psi)|^2_{1,K,EI},\qquad\forall\,w,\psi\in V^1_{K,EI,0} (0,\ell ),
\end{equation}
where
\begin{align}\label{eq:cost_pcr}
C_{D,K,EI}:=
  \ell^2 \max\Big\{ \frac{1}{K(\ell)},\; \frac{1}{EI(\ell)}\Big\}
\min\Big\{4,\frac{1}{2-\mu}\Big\}
 \max\Big\{2,\;
 1+\frac{2\ell^2 \|K\|_\infty}{K(\ell)}\min\Big\{4,\frac{1}{2-\mu}\Big\}
\Big\}.
 \end{align}

\end{proposition}
\begin{Proof}
By applying Proposition 2.2 from \cite{AlabauCannarsaLeugering2017} with adjusted constants and \eqref{eq:a}, we individually have
\begin{align}
&\int\limits_0^\ell w^2(x)\mathrm dx\leq
C_{D1}\int_0^\ell K(x)w'(x)^2 \mathrm  dx \leq
2C_{D1}\int\limits_0^\ell K(x)(w'(x)+\psi(x))^2\mathrm  dx +2C_{D1}\|K\|_\infty\int\limits_0^\ell\psi^2(x) \mathrm dx
\end{align}
and
\begin{align}\int\limits_0^\ell \psi^2(x) \mathrm  dx\leq C_{D2}\int_0^\ell EI(x)\psi'(x)^2 \mathrm  dx,\end{align}
where
\begin{align*}
& C_{D1} := \frac{\ell^2}{K(\ell)}\min\Big\{4,\frac{1}{2-\mu_{K}}\Big\},
\\
& C_{D2} := \frac{\ell^2}{EI(\ell)}\min\Big\{4,\frac{1}{2-\mu_{EI}}\Big\}.
\end{align*}
Putting things together, we arrive at
$\|(w,\psi)\|^2_{L^2(0,\ell)^2}\leq \tilde{C}_{D,K,EI}\; |(w,\psi)|^2_{1,K,EI}$, where
\begin{align*}\label{eq:cost_pcr}
\tilde{C}_{D,K,EI}:= &\;  \ell^2 \max\Big\{ \frac{1}{K(\ell)},\; \frac{1}{EI(\ell)}\Big\} \max\Big\{2
\min\Big\{4,\frac{1}{2-\mu_K}\Big\}, \notag
\\  &  \qquad   \qquad
\Big(1+\frac{2\ell^2 \|K\|_\infty}{K(\ell)}\min\Big\{4,\frac{1}{2-\mu_K}\Big\}\Big)  \min\Big\{4,\frac{1}{2-\mu_{EI}}\Big\}
\Big\}.
 \end{align*}
 It is clear that $C_{D,K,EI}\le \tilde{C}_{D,K,EI}$. The proof is completed.
\end{Proof}

\begin{example}\label{exa:power}\rm The following example is prototypical for \eqref{eq:hp_a}.
Let  $\theta\in[0,2)$ be given. Define
\begin{equation}\label{eq:power}
a(x)=x^\theta,\;\; \forall x\in[0,\ell ].
\end{equation}
In this case, we have $\|u\|^2_{L^2(0,\ell )}\leqslant\ell^{2-\theta} \min \Big\{ 4, \dfrac{1}{ 2-\theta}\Big\}\; |u|^2_{1,a},\;\forall u\in V^1_{a,0} (0,\ell )$. Hence, if $K(x)=EI(x)=x^\theta$ and $\ell=1$, we have  $\theta\in [0,\frac{7}{4}),$  $C_{D,K,EI}=\frac{4-\theta}{(2-\theta)^2}$,
 while
 $C_{D,K,EI}= 36 $  for  $\theta\in [\frac{7}{4},2)$,
 else.
\end{example}

Next, we define
$$
V^2_{K,EI}(0,\ell )=\big\{(w,\psi) \in V^1_{K,EI}(0,\ell )\; \big|\; ~K(w'+\psi), EI \psi' \in H^1(0,\ell )\big\},
$$
where $H^1(0,\ell )$ denotes the classical Sobolev space of all functions $u\in L^2(0,\ell )$ such that $u'\in L^2(0,\ell )$. Notice that, if $(w,\psi)\in V^2_{K,EI}(0,\ell )$, then $Kw', EI\psi'$ are continuous on $[0,\ell ]$.

\begin{remark}\label{Poincare_Neumann}
Clearly, in Proposition \ref{Poincare_inequality} we can also discriminate between the cases where $K$ or $EI$ exhibits strong degeneracy. Notice that in the case of weak degeneracy, the corresponding $a=K$ and or $a=EI$ satisfies $ \frac{^1}{a(x)}\in L^1(0,\ell)$ and this makes it possible to use the same proof of the Poincar\'{e} inequality, even in the case of Neumann controls.
If, however,  we take Neumann controls and consider strong degeneration for, say, $K$, then, there is no Dirichlet condition at $x=0$ to be used in the proof of the corresponding Poincar\'{e} inequality.
In that case, we need to resort to a result by Chua and Wheeden  \cite{ChuaWheeden00}
 and restrict ourselves to functions satisfying $\int\limits_0^\ell f(x)dx =0$ and we introduce the quotient space
$$
\tilde{V}^2_{K,EI}(0,\ell)=\Big\{u \in V^1_{K,EI}(0,\ell )~\Big|~K(w'+\psi), EI\psi' \in H^1(0,\ell), \int\limits_0^\ell w(x) dx =0, \, \int\limits_0^\ell \psi(x) dx =0 \Big\}.
$$
\end{remark}
Another way to handle the strong degenerate case is to introduce mixed boundary conditions or Robin-type boundary conditions at $x=\ell$. This will turn out to be important anyway when using state and velocity boundary feedback controls as introduced later. In this case, we imposes the boundary controls as in \eqref{Timo_weak_strong_Neumann}
and introduce the extended bilinear form
\begin{align}\label{bilinear_gd}
    B_{\gamma,\delta}((z,\xi),(u,\chi))=& \int\limits_0^\ell\left\{K(x)(z_x(x)+\xi(x))(u_x(x)+\chi(x))+EI(x)\xi_x(x)\chi_x(x)
    \right\}dx
   \notag \\ & + \gamma z(\ell)u(\ell)+\delta \xi(\ell)\chi(\ell).
\end{align}
This gives rise to the norm in $V^1_{K,EI}(0,l):$ $$ \|| (w,\psi)|\|_{\gamma,\delta}:=\sqrt{B_{\gamma,\delta}((w,\psi),(w,\psi))}, \quad \forall\;  (w,\psi)\in V^1_{K,EI}(0,l).$$
\begin{proposition}\label{Poincare_gd}
   For $(\w,\psi)\in V_{K,EI}^1(0,\ell)$ and $K,{EI}\in \mathcal{A}$, we have
   \begin{align}\label{Poincare_inequaltity_gd}
     \|(w,\psi)\|_{L^2(0,\ell)^2}^2\leq C_{N,K, EI}\||(w,\psi)\||^2_{\gamma,\delta},
   \end{align}
 where $C_{N,K,EI}$ is given by
   \begin{align}\label{Poincare_Robin}
      & \tilde{C}_{N1}:=   2\ell^2 \max\Big\{ \frac{1}{K(\ell)},\; \frac{1}{EI(\ell)}\Big\} \min\Big\{2,\frac{1}{2-\mu}\Big\} \max\Big\{2,\;
\Big(1+\frac{2\ell^2 \|K\|_\infty}{K(\ell)}\min\Big\{2,\frac{1}{2-\mu}\Big\}\Big) \Big\},\notag\\
      & \tilde{C}_{N2}:= \max\Big\{ \frac{2\ell}{\gamma}, \; \frac{8\ell^3\|K\|_\infty}{\delta K(\ell)}\min\Big\{2,\frac{1}{2-\mu}\Big\} \Big\}, \;\;  \notag
      \\
      &C_{N,K,EI}:= \max\big\{\tilde{C}_{N1},\; \tilde{C}_{N2}\big\}.\end{align}

\end{proposition}
\begin{Proof}
The proof is analogous to the one of the Proposition \ref{Poincare_inequality}. The difference is that we now have to account for the values of $w,\psi$ at $x=\ell$. We obtain,
\begin{align*}
    \int\limits_0^\ell w^2(x)\mathrm dx \leq& C_{N1}\int_0^\ell K(x)w'(x)^2 \mathrm  dx+2\ell w^2(\ell)\\
    \leq &2 C_{N1}\int\limits_0^\ell K(x)(w'+\psi)^2\mathrm dx+2\ell w^2(\ell)  +2C_{N1}\|K\|_{\infty} \int_0^\ell  \psi^2(x) \mathrm  dx ,
\end{align*}
and
\begin{align*}
    \int\limits_0^\ell \psi^2(x) \mathrm dx\leq C_{N2} \int_0^\ell EI(x)\psi'(x)^2 \mathrm  dx+2\ell \psi^2(\ell),
\end{align*}
where
\begin{align*}
&C_{N1} :=\frac{2\ell^2}{K(\ell)}\min\Big\{\frac{1}{2-\mu_{K}},2\Big\},
\\ &
C_{N2} :=\frac{2\ell^2}{EI(\ell)}\min\Big\{\frac{1}{2-\mu_{EI}},2\Big\}.
\end{align*}
Putting things together, we arrive at $ \|(w,\psi\|_{L^2(0,\ell)^2}^2\leq \tilde{C}_{N,K, EI}\||(w,\psi)\||^2_{\gamma,\delta}$ where
\begin{align*}
\tilde{C}_{N,K, EI}:= \max\Big\{2C_{N1},\;  (1+2C_{N1}\|K\|_{\infty})C_{N2},\;
{{2\ell}\over{\gamma}},\;  (1+2C_{N1}\|K\|_{\infty}) {{2\ell}\over{\delta}}\Big\}.
 \end{align*}
\end{Proof}

The well-posedness of the homogeneous system \eqref{Timo_weak_strong_Neumann} with boundary conditions \eqref{Timo_weak_strong_Neumann_homogen} or \eqref{Timo_weak_strong_Dirichlet_homogen} instead of \eqref{Timo_weak_strong_Neumann}$_5$ follows by standard arguments outlined e.g. in \cite{AlabauCannarsaLeugering2017}. We, therefore, just give an outline here and leave the details to the reader. In the case of boundary controls, we use the method of transposition.

 \subsection{Well-posendess}\label{well-posendess}
\subsubsection{homogeneous Dirichlet case}
 We introduce the space $(\mathcal{H},\|\cdot \|)$ such that
 $$
 \mathcal{H}:=  L^2(0,\ell)^2  , \; \;\|(v,\phi)\|^2_{\mathcal {H}}=\int\limits_0^\ell \left( \rho v^2+I_\rho \phi^2\right)dx.
 $$
 Moreover, for weak degeneracy, we define $(\mathcal{V}_{D,w}, \|\cdot \|_{\mathcal{V}_{D,w}})$ as
 $$
\mathcal{V}_{D,w}:=\{(v,\phi)\in V^1_{K,EI,0}(0,\ell)\; |\;  v(0)=\phi(0)=0\}\;\;
 \mbox{with} \;\;
  \|(v,\phi)\|_{V_{D,w}}:=(B(v,\phi),(v,\phi))^{1\over2},$$
 while for strong degeneracy we set
 $$
\mathcal{V}_{D,s}:=V^1_{K,EI,0}(0,\ell)\;\; \mbox{with} \;\;
 \|(v,\phi)\|_{V_{D,s}}:=(B(v,\phi),(v,\phi))^{1\over2},$$
 where in both expressions the bilinear form $B$ is given by \eqref{bilinear_form}.

 We define the spatial operator as
 \begin{equation}
 \label{defa}
 A(v,\phi):= \left((K(v_x+\phi))_x, (EI \phi_x)_x-K(v_x+\phi)\right).
 \end{equation}
 In the case of weak degeneracy, $A$ acts on
$$
 D(A)_{D,w}:=\{(v,\phi)\in \mathcal{V}_{D,w}\,|\, A(w,\phi)\in L^2(0,\ell)^2 \},
 $$
 while in the strongly degenerated case it acts on
 $$
 D(A)_{D,s}:=\{(v,\phi)\in \mathcal{V}_{D,s}\;|\; A(v,\phi)\in L^2(0,\ell)^2,
 \lim\limits_{x\rightarrow 0^+}K(x)(w_x(x)+\psi(x))=0,\;  \lim\limits_{x\rightarrow 0^+}EI(x)\psi_x(x)=0
 \}.
 $$

\noindent
Then, $
\left(A(v,\phi),(w,\psi)\right)=B\left((v,\phi),(w,\psi)\right),
$
for all $(v,\phi), (w,\psi)$ in the corresponding spaces $D(A)_{D,w}$ or $D(A)_{D,s}$.

In both cases it is easily to prove that \( A \) is the generator of a $C_0$ contraction semigroup   in \( \mathcal{V}_{D,s} \times \mathcal{H}  \) or
\( \mathcal{V}_{D,w} \times \mathcal{H}  \). The system \eqref{Timo_weak_strong_Neumann}$_{1,2,3,4}$ with \eqref{Timo_weak_strong_Dirichlet_homogen} can be written as
\begin{equation}
\label{ab}
\begin{cases}
U'(t) = AU(t),\;\; t \geq 0, \\
U(0) = U_0.
\end{cases}
\end{equation}
For the initial conditions  \( (w^0, \psi^0, w^1,\psi^1 ) \in \mathcal{V}_{D,w/s} \times \mathcal{H} \), there exists unique mild solution
to \eqref{ab} such that
$$
(u,\psi) \in C^1([0, \infty); {\cal H}) \cap C([0, \infty); \mathcal{V}_{D,w})
$$
Furthermore, when the initial conditions are \( (w^0, \psi^0, w^1,\psi^1 ) \in (V^2_{K,EI} (0,\ell) \cap \mathcal{V}_{D,w/s}) \times \mathcal{V}_{D,w/s}  \), a unique classical solution to equation \eqref{ab} exists, with the following properties:
\[
(u,\psi) \in C^2([0, \infty); {\cal H}) \cap C^1([0, \infty);\mathcal{V}_{D,w}) \cap C([0, \infty); V^2_{K,EI} (0,\ell) \cap \mathcal{V}_{D,w/s} ).
\]


 \subsubsection{The homogeneous Neumann case}\label{section_NB}
 In this case, for weak degeneracy, we define $(\mathcal{V}_{N,w}, \|\cdot \|_{\mathcal{V}_{N,w}})$ as
 $$
\mathcal{V}_{N,w}:=\{(v,\phi)\in V^1_{K,EI}(0,\ell)\,|\, v(0)=\phi(0)=0\} \;\; \mbox{ with} \;\;
  \|(v,\phi)\|_{V_{N,w}}:=(B(v,\phi),(v,\phi))^{1\over2},
 $$
 while for strong degeneracy we work in the quotient space and set
 $$
\mathcal{V}_{N,s}:=\Big\{(v,\phi)\in V^1_{K,EI}(0,\ell)\,\Big|\, \int\limits_0^\ell vdx=\int\limits_0^\ell \phi dx=0\Big\}
\;\; \mbox{with}\;\;
 \|(v,\phi)\|_{V_{D,s}}:=(B(v,\phi),(v,\phi))^{1\over2}.$$

 The action of $A$ is the same as \eqref{defa} while its domains are now given as follows.
 In the case of weak degeneracy, $A$ acts on
 $$
 D(A)_{N,w}:=\{(v,\phi)\in \mathcal{V}_{N,w}| A(v,\phi)\in L^2(0,\ell)^2\},
 $$
 while in the strongly degenerated case it acts on
 \begin{align*}
 D(A)_{N,s}:=\{(v,\phi)\in \mathcal{V}_{N,s}\; | \;& A(v,\phi)\in L^2(0,\ell)^2,
 \lim\limits_{x\rightarrow 0^+}K(x)(w_x(x)+\psi(x))=0,\;  \lim\limits_{x\rightarrow 0^+}EI(x)\psi_x (x)=0,\\&
 K(\ell)(w_x(\ell)+\psi(\ell)) =0,\;  EI(\ell)\psi_x (\ell)=0
 \}.
 \end{align*}
 Again,  in both cases it is easily seen that $A$ generates strongly continuous semigroup in $\mathcal{V}_{N,w/s}\times \mathcal{H}$ such that the solutions to \eqref{Timo_weak_strong_Neumann}$_{1,2,3,4}$ with \eqref{Timo_weak_strong_Neumann_homogen} and $\gamma=\delta=0$ are well-defined and can be classified into classical, strong, mild and weak solutions.

 \subsubsection{The Robin case}
  In this case we have $\gamma,\delta>0$. For weak degeneracy, we define $(\mathcal{V}_{R,w}, \|\cdot\|_{\mathcal{V}_{R,w}})$ as
  $$
\mathcal{V}_{R,w}:= \mathcal{V}_{N,w}  \;\; \mbox{with}
\;\;   \|(v,\phi)\|_{V_{R,w}}:=B_{\gamma,\delta}((v,\phi),(v,\phi))^{1\over2},
$$
 while for strong degeneracy, we set
$$
\mathcal{V}_{R,s}:=V^1_{K,EI}(0,\ell),
\;\;  \mbox{with} \|(v,\phi)\|_{V_{R,s}}:=B_{\gamma,\delta}((v,\phi),(v,\phi))^{1\over2}.
 $$

The action of $A$ is the same as \eqref{defa} while its domains are now given as follows.
 In the case of weak degeneracy, $A$ acts on
 $$
 D(A)_{R,w}:=\{(v,\phi)\in \mathcal{V}_{R,w}| A(v,\phi)\in L^2(0,\ell)^2\},
 $$
 while in the strongly degenerated case it acts on
\begin{align*}
 D(A)_{R,s}:=\{(v,\phi)\in \mathcal{V}_{R,s}\;|\;& A(v,\phi)\in L^2(0,\ell)^2,
 \lim\limits_{x\rightarrow 0^+}K(x)(w_x(x)+\psi(x))=0,\;  \lim\limits_{x\rightarrow 0^+}EI(x)\psi_x (x)=0,\\&
 K(\ell)(w_x(\ell)+\psi(\ell))+\gamma w(\ell) =0,\;  EI(\ell)\psi_x (\ell)+\delta \psi(\ell)=0\}.
 \end{align*}
 Again,  in both cases it is easily seen that $\mathcal{A}$ generates strongly continuous semigroup in $\mathcal{V}_{R,w/s}\times \mathcal{H}$ such that the solutions to \eqref{Timo_weak_strong_Neumann}$_{1,2,3,4}$ with \eqref{Timo_weak_strong_Neumann_homogen} are well-defined and can be classified into classical, strong, mild and weak solutions.


\section{Multiplier identities}\label{multiplier_identities}

\subsection{Multiplier identities}
Before proving the main results of this paper, we first present several multiplier equalities, which lead to some direct and indirect inequalities.

\begin{proposition}\label{main_multiplier_identiy}
Let $K,EI\in \mathcal{A}$ and let $(w,\psi)$ be a pair of smooth functions satisfying the Timoshenko beam equations \eqref{Timo_weak_strong_Neumann}$_{1,2}$
Then, for $S\leq T$ the following identity holds.
\begin{align}\label{multiplier_identity_main}
& \int\limits_S^T\ell  \left[ \rho w_t^2( \ell)+I_\rho \psi_t^2( \ell) \right] dt+ \int\limits_S^T \ell \left[K(\ell)(w_x( \ell)+\psi( \ell))^2+EI(\ell) \psi_x^2( \ell) \right]dt
\notag \\
 = \;& 2(F(T)-F(S))+\int\limits_S^T\int\limits_0^\ell \left[ \rho w_t^2+I_\rho \psi_t^2 + (K-xK')(w_x+\psi)^2+(EI-xEI')\psi_x^2\right]dx dt
 \\
& \quad\quad +2\int\limits_S^T\int\limits_0^\ell x\left[-\rho  w_t \psi_t + K(w_x+\psi)\psi_x \right]dx dt,\notag
\end{align}
where
$$
F(t) :=\int\limits_0^\ell x\left[\rho w_t (w_x +\psi )+I_\rho \psi_t \psi_x \right]dx
$$
\end{proposition}
\begin{Proof}
We multiply the first equation in \eqref {Timo_weak_strong_Neumann} by $xw_x$ and the second by $x\psi_x$ and integrate by parts to get
\begin{align}\label{multipliers1}
0=\;& \int\limits_S^T\int\limits_0^\ell \left[\left(\rho w_{tt}-(K(w_x+\psi))_x\right)) xw_x+\left(I_\rho\psi_{tt}-(EI\psi_x)_x+K(w_x+\psi)\right) x \psi_x\right] dx dt\notag\\
=\;& \int\limits_0^\ell\left( \rho x w_x w_t+I_\rho x \psi_x\psi_t \right)|_S^Tdx-\int\limits_S^T\int\limits_0^\ell \left[x\rho w_t w_{xt}+xI_\rho \psi_t\psi_{xt}\right] dx dt\notag\\
& -\int\limits_S^T\int\limits_0^\ell x\left[(K(w_x+\psi))_x w_x+(EI\psi_x)_x \psi_x-K(w_x+\psi)\psi_x\right] dx dt\notag\\
=\;& F_1(T)-F_1(S)-\frac{1}{2}\int\limits_S^T\int\limits_0^\ell x \left[ \rho(w_t)^2_x+I_\rho(\psi_t)^2_x\right] dx dt+I.
\end{align}
where
$$
F_1 := \int\limits_0^\ell\left( \rho x w_x w_t+I_\rho x \psi_x\psi_t \right) dx,
$$
and
$$
I := -\int\limits_S^T\int\limits_0^\ell x \left[\big(K(w_x+\psi)\big)_x w_x+(EI\psi_x)_x \psi_x-K(w_x+\psi)\psi_x\right] dx dt.
$$
 Note that
 $$
 \int\limits_S^T\int\limits_0^\ell  x (K (w_x+\psi))_x\psi dx dt  = \int\limits_0^\ell \rho x  w_t \psi|_S^Tdx-\int\limits_S^T\int\limits_0^\ell \rho x w_t \psi_t dx  .
 $$
We have, after some elementary calculations and integration by parts,
\begin{align}\label{multipliers2}
I
=&  -\frac{1}{2}\int\limits_S^T x\left(K(w_x+\psi)^2+EI \psi_x^2\right)|_0^\ell dt\notag \\
& +\frac{1}{2}\int\limits_S^T\int\limits_0^\ell  \left[(K-xK')(w_x+\psi)^2+(EI-xEI')\psi_x^2\right] dx dt\notag\\
&  +\int\limits_0^\ell x \rho w_t \psi|_S^Tdx-\int\limits_S^T\int\limits_0^\ell \rho x w_t \psi_t dx dt+\int\limits_S^T\int\limits_0^\ell xK(w_x+\psi)\psi_x dx dt.
\end{align}
Define $ F_2(t):= \int\limits_0^\ell x\rho w_t(t,x)\psi(t,x)dx.$
Then
\begin{align}\label{F}
F(t) =  F_1(t)+F_2(t).
\end{align}
We go back to \eqref{multipliers1} and obtain
\begin{align*}
0 =& F_1(T)-F_1(S)-\frac{1}{2}\int\limits_S^T \ell  \left[ \rho w_t^2(\ell)+I_\rho\psi_t^2(\ell) \right] dt +\frac{1}{2}\int\limits_S^T\int\limits_0^\ell  \left[ \rho w_t^2+I_\rho \psi_t^2\right] dx dt+I.
\end{align*}
Therefore, \eqref{multipliers1},\eqref{multipliers2} and \eqref{F} imply
\begin{align*}
& \frac{1}{2}\int\limits_S^T\ell  \left[ \rho w_t^2( \ell)+I_\rho \psi_t^2( \ell) \right] dt+ \frac{1}{2}\int\limits_S^T \ell \left[K(\ell)(w_x( \ell)+\psi( \ell))^2+EI(\ell) \psi_x^2( \ell) \right]dt\notag\\
& = F(T)-F(S)+\frac{1}{2}\int\limits_S^T\int\limits_0^\ell \left[ \rho(w_t)^2+I_\rho(\psi_t)^2 + (K-xK')(w_x+\psi)^2+(EI-xEI')\psi_x^2\right]dx dt\\
& +\int\limits_S^T\int\limits_0^\ell x\left[-\rho  w_t \psi_t + K(w_x+\psi)\psi_x \right]dx dt,
\end{align*}
as desired.
\end{Proof}
We introduce the total energy $E(t)$ at time $t$ for system \eqref{Timo_weak_strong_Neumann}$_{1,2,3,4}$ with \eqref{Timo_weak_strong_Dirichlet_homogen}:
\begin{equation}\label{energy}
E(t):=\frac{1}{2}\int\limits_0^\ell \left[\rho w_t^2(t,x)+I_\rho\psi(t,x)^2+K(w_x(t,x)+\psi(t,x))^2+EI\psi_x^2(t,x) \right]dx ,
\end{equation}
and for \eqref{Timo_weak_strong_Neumann}$_{1,2,3,4}$ with \eqref{Timo_weak_strong_Neumann_homogen}:
\begin{equation}\label{E-gd}
E_{
\gamma,\delta
}(t):=\frac{1}{2}\Big[\int\limits_0^\ell \left[\rho w_t^2(t,x)+I_\rho\psi(t,x)^2+K(w_x(t,x)+\psi(t,x))^2+EI\psi_x^2(t,x) \right]dx
+ \gamma w^2(t,\ell)+\delta \psi^2(t,\ell)\Big],
\end{equation}
Taking the time derivative and integration by parts implies
\begin{align}\label{dEdt}
\frac{d}{dt}E(t)=K(x)(w_x(t,x)+\psi(t,x))w_t(t,x)|_0^\ell+ EI(x)\psi_x(t,x)\psi_t(t,x)|_0^\ell.
\end{align}
\begin{proposition}\label{conservation}
Let $(w,\psi)$ be a pair of functions satisfying \eqref {Timo_weak_strong_Neumann}$_{1,2}$ in the classical sense, then, the energy identity \eqref{dEdt} holds. If in addition $(w,\psi)$ satisfies \eqref{Timo_weak_strong_Dirichlet_homogen} or \eqref{Timo_weak_strong_Neumann_homogen}, then the energy is conserved. In particular,
\[
E(t)=E(0), \; \; \mbox{  or }\;\;  E_{\gamma,\delta}(t)=E_{\gamma,\delta}(0) \; \;\; \forall\;  t\in (0,T).
\]
\end{proposition}
We now establish a second multiplier identity.
\begin{proposition}\label{multiplier_identity2}
Let $(w,\psi)$ be a pair of functions satisfying \eqref {Timo_weak_strong_Neumann}$_{1,2}$. Then, for $S\leq T$, the following identity holds.
\begin{align}\label{multiplier_identity3}
0=&\int\limits_S^T\int\limits_0^\ell \left[ K(w_x+\psi)^2+EI\psi_x^2- \left(\rho w_{t}^2+I_\rho \psi_t^2\right)\right] dx dt +\int\limits_0^\ell \left( \rho w_t w +I_\rho \psi_t  \psi \right) dx|_S^T\\
& \quad\quad\quad -\int\limits_S^TK (w_x +\psi )w +EI(x)\psi_x \psi dt|_0 ^\ell.\notag
\end{align}
If moreover $(w,\psi)$ satisfies \eqref{Timo_weak_strong_Dirichlet_homogen} or \eqref{Timo_weak_strong_Neumann_homogen}, then we have
\begin{align}\label{multiplier_identity4}
0=&\int\limits_S^T\int\limits_0^\ell \left[ K(w_x+\psi)^2+EI\psi_x^2- \left(\rho w_{t}^2+I_\rho \psi_t^2\right)\right] dx dt+\int\limits_0^\ell \left( \rho w_t w +I_\rho \psi_t  \psi \right) dx|_S^T,
\end{align}
or \textcolor{blue}{
\begin{align}\label{multiplier_identity5}
0=&\int\limits_S^T\int\limits_0^\ell \left[ K(w_x+\psi)^2+EI\psi_x^2- \left(\rho w_{t}^2+I_\rho \psi_t^2\right)\right] dx dt+\int\limits_0^\ell \left( \rho w_t w +I_\rho \psi_t  \psi \right) dx|_S^T \notag
\\ &+ \int\limits_S^T \left[\gamma w^2(\ell) + \delta \psi^2(\ell)\right]dt.
\end{align}}
\end{proposition}
\begin{Proof}
We multiply \eqref {Timo_weak_strong_Neumann}$_{1,2}$ by $w$ and $\psi$, respectively. The calculations are standard.
\end{Proof}

\subsection{The direct and indirect inequalities}\label{inequalities}
 We now embark on estimates for the traces of $w_t,w_x,\psi_t,\psi_x$ at $(t,\ell)$ using the multiplier identities \eqref{multiplier_identity_main}, \eqref{multiplier_identity4} and \eqref{multiplier_identity5}. To this end, we finally have to focus on one of the boundary condition, namely Dirichlet, Robin or Neumann conditions.  We first consider the term $F(t)$.
 \begin{proposition}\label{F_estimate}
 Let $(w,\psi)$ be a pair of functions satisfying \eqref {Timo_weak_strong_Neumann}$_{1,2}$,
 $\mu = \max\{\mu_K,\; \mu_{EI}\} $ and $F$ be given by \eqref{F}. We consider the case of weak degeneracy, where in this case we may include $\mu_K, \mu_{EI}=1$, and strong degeneracy separately.
 \begin{enumerate}
     \item $\mu_K, \mu_{EI}\in [0,1]:$
\begin{equation}\label{F_estimate_1}
 F(t)\leq \ell \max\Big\{ \sqrt{\frac{\rho}{K(\ell)}}, \sqrt{\frac{I_\rho}{EI(\ell)}}\Big\}E(t).
 \end{equation}
 \item $\mu_K,\mu_{EI}\in (1,2):$
 \begin{equation}\label{F_estimate_2}
 F(t)\leq \begin{cases}
  \displaystyle\ell^{2-\mu} \max\Big\{ \sqrt{\frac{\rho}{K(\ell)}}, \; \sqrt{\frac{I_\rho}{EI(\ell)}}\Big\}E(t) & \text{ if } \ell \leq 1
,\\
   \displaystyle
    \ell^{\mu} \max\Big\{ \sqrt{\frac{\rho}{K(\ell)}},  \; \sqrt{\frac{I_\rho}{EI(\ell)}}\Big\}E(t) &  \text{ if } \ell \geq 1.\end{cases}
 \end{equation} \end{enumerate}

 We denote by $C_F$
 the corresponding constants in \eqref{F_estimate_1} and \eqref{F_estimate_2}.

 \begin{Proof}
In the first case, we have
\begin{align*}
   & \int\limits_0^\ell x\rho w_t(t,x)(w_x(t,x)+\psi(t,x))dx=
   \sqrt{\frac{\rho\ell}{K(\ell)}}
   \int\limits_0^\ell   \sqrt{x}\sqrt{\rho}w_t(t,x)\sqrt{x K(\ell)\over\ell}(w_x(t,x)+\psi(t,x) )dx\\
   & \qquad\qquad \leq \frac{1}{2}\ell \sqrt{\frac{\rho}{K(\ell)}}\int\limits_0^\ell[ \rho w_t^2(t,x)+\frac{x}{\ell}K(\ell)(w_x(t,x)+\psi(t,x))^2]dx\\
   &\qquad\qquad \leq \frac{1}{2}\ell \sqrt{\frac{\rho}{K(\ell)}}\int\limits_0^\ell [\rho w_t^2(t,x)+K(x)(w_x(t,x)+\psi(t,x))^2]dx.
\end{align*}
Similarly,
\begin{align*}
  \int\limits_0^\ell xI_\rho \psi_t(t,x)\psi_x(t,x) dx\leq\frac{1}{2}\ell \sqrt{\frac{I_\rho}{EI(\ell)}}\int\limits_0^\ell [I_\rho \psi_t^2(t,x)+EI(x) \psi_x^2(t,x)]dx.
\end{align*}
 Thus, putting these inequalities together, we obtain \eqref{F_estimate_1}.

For the second case, we have to use the inequalities  $ (\frac{x}{\ell})^{\mu_K}K(\ell)\leq K(x)$  and, accordingly, $(\frac{x}{\ell})^{\mu_{EI}}EI(\ell)\leq EI(x)$. Then,
\begin{align*}
 & \int\limits_0^\ell x\left\{\rho w_t(t,x)(w_x(t,x)+\psi(t,x))\right\} dx\notag\\
 =& \sqrt{\frac{\rho}{K(\ell)}}\int\limits_0^\ell
 \sqrt{\rho}x^{\frac{2-\mu_K}{2}}w_t(t,x)x^{\frac{\mu_K}{2}}\sqrt{K(\ell)}(w_x(t,x)+\psi(t,x))dx \\
\leq &  \frac{1}{2}\sqrt{\frac{\rho}{K(\ell)}}\int\limits_0^\ell [\rho \ell^{2-\mu_k}w_t^2(t,x)+\ell^{\mu_k}K(x)(w_x(t,x)+\psi(t,x))^2]dx\\
\leq & \frac{1}{2}\sqrt{\frac{\rho}{K(\ell)}}\max\{\ell^{2-\mu_K},\ell^{\mu_K}\}
\int\limits_0^\ell[ \rho w_t^2(t,x)+K(x)(w_x(t,x)+\psi(t,x))^2]dx.
 \end{align*}
 Similarly,
 \begin{align*}
     \int\limits_0^\ell x I_\rho \psi_t(t,x)\psi_x(t,x)dx\leq \frac{1}{2}\sqrt{\frac{\rho}{K(\ell)}}
     \max\{\ell^{2-\mu_{EI}},\ell^{\mu_{EI}}\}
     \int\limits_0^\ell [ I_\rho \psi_t^2(t,x)+EI \psi_x^2(t,x)]dx. \end{align*}
     Hence, depending on $\ell \leq 1$ or $\ell \geq 1$, we arrive at \eqref{F_estimate_2}.
 \end{Proof}
 \end{proposition}

From now on until Section \ref{stabilizability}, we take $S=0$. As for the last term in \eqref{multiplier_identity_main}, we have
 \begin{align}\label{final_term}
 & \int\limits_0^T\int\limits_0^\ell x\left(-\rho  w_t \psi_t + K(w_x+\psi)\psi_x \right)dx dt\notag\\
 = \quad &\int\limits_0^T\int\limits_0^\ell   \left(
  -x\sqrt{\frac{\rho}{ I_\rho }}  w_t \sqrt{\rho I_\rho}\psi_t +\sqrt{\frac{K(x)  }{EI(\ell)}}(w_x+\psi)\sqrt{K(x)}\sqrt{EI(\ell) } x \psi_x \right)dx dt
  \end{align}
  \begin{align*}
  & \leq  \frac{1}{2}\int\limits_0^T\int\limits_0^\ell\left\{  \ell
  \sqrt{\frac{\rho}{ I_\rho }}\left(\rho w_t^2+ I_\rho \psi_t^2\right)+
  \sqrt{\frac{\|K\|_\infty }{EI(\ell)}}\left(K(x)(w_x+\psi)^2+\ell^2 EI(x)\psi_x^2\right)\right\} dx dt\notag
   \\&
  \leq C_h \int_0^TE(t)dt.
  \notag
 \end{align*}
 where we use $EI(\ell)x^2 \le \ell^2 EI(x)$ in Proposition \ref{re:a}, and
 $$
 C_h:=\max\Big\{
  \ell
  \sqrt{\frac{\rho}{ I_\rho }},\;  \sqrt{\frac{\|K\|_\infty }{EI(\ell)}}
  \max\big\{1,\ell^2)\big\}\Big\}.
 $$
 \begin{remark}\label{smallh}
     We recall the definition of $K(x)=G(x)h$ and $\rho=m_0 h$. Hence $\frac{K(x)}{EL(\ell)}=\frac{G(x)h}{EI(\ell)}\leq h\frac{\|G\|_\infty }{EI(\ell)}$ and this can be made sufficiently small for small $h$. Similarly, $\frac{\rho}{I_\rho}=\frac{m_0 h}{m_0 I}=\frac{h}{I}$ can be made sufficiently small for small $h$. Thus $C_h$
     can be made sufficiently small for small $h$.
     \end{remark}

We are now in the position to prove upper bounds for the traces involved.
\begin{proposition}\label{upper_Dirichlet} Assume   $\mu_K, \; \mu_{EI} \in [0,2)$.
\begin{enumerate}
\item[{\rm (i)}]
 Let $(w,\psi)$ be classical solutions of \eqref{Timo_weak_strong_Neumann} with homogeneous Dirichlet conditions \eqref{Timo_weak_strong_Dirichlet_homogen}. Then,
\begin{equation}\label{upper_Dirichlet_estimate}
\int\limits_0^T \ell \left(K(\ell)w_x^2(t,\ell)+ EI(\ell) \psi^2_x(t,\ell)\right)dt\leq C_{BC}E(0),
\end{equation}
where
$$
C_{BC} :=4C_F   +2(1+\mu) + 2C_hT.
$$

\item[{\rm (ii)}]  Let $(w,\psi)$ be the classical solution of \eqref{Timo_weak_strong_Neumann} with homogeneous Robin conditions \eqref{Timo_weak_strong_Neumann_homogen}. Then,
\begin{equation}\label{upper_Neumann_estimate}
\int\limits_0^T \ell \left(\rho w_t^2(t,\ell)+ I_\rho  \psi_t^2(t,\ell) + \frac{\gamma^2}{K(\ell)}  w^2(t,\ell) + \frac{\delta^2}{EL(\ell)} \psi^2(t,\ell)\right)dt\leq C_{BC}E(0).
\end{equation}

\end{enumerate}
\end{proposition}
\begin{Proof}
We use the multiplier identity \eqref{multiplier_identity_main} (with $S=0)$ and the estimates \eqref{F_estimate_1}, \eqref{F_estimate_2}, \eqref{final_term},
\begin{align*}
&\int\limits_0^T \ell \left(K(\ell)w_x^2(t,\ell)+ EI(\ell) \psi^2_x(t,\ell)\right)dt\leq
4C_F E(0) +2T(1+\mu)E(0)+ 2C_hT E(0)
\end{align*}
The estimation is analogous for the Neumann case.
\end{Proof}
\begin{remark}
 At this point, we have to discriminate between four cases.
 \begin{enumerate}
 \item $\mu_K, \mu_{EI} \in [0,1)$ (weak degeneracy) ,
 \item $\mu_K\in [0,1), \mu_{EI}\in [1,2)$ (mixed degeneracy: weak-strong),
 \item $\mu_K\in [1,2), \mu_{EI}\in [0,1)$ (mixed degeneracy: strong-weak),
 \item $\mu_K\in [1,2), \mu_{EI}\in [1,2)$ (strong degeneracy).
 \end{enumerate}
 Then in case 1., we can take Dirichlet or Neumann conditions at $x=\ell$,as in this case $\frac{1}{K},\frac{1}{EI}$ are both integrable. In all other cases, we have to work in a quotient space if we apply Neumann controls, while Dirichlet controls at $x=\ell$ can be used without resorting to quotient spaces. In order to keep matters simple, we focus on cases 1. and 4. only, where in case 4., we apply Dirichlet controls. We remark already at this point that the term \eqref{final_term} makes it impossible to allow for the full range of degeneracy, i..e $\mu_K, \mu_{EI}$ very close to 2. For that matter, weak degeneracy is more realistic anyway.
 \end{remark}

 The next estimates concern the last term in \eqref{multiplier_identity4}.
 \begin{proposition}\label{lower_order_estimate}
Let the case 1. together with \eqref{Timo_weak_strong_Dirichlet_homogen}  or \eqref{Timo_weak_strong_Neumann_homogen} hold or cases 2.-4. along with \eqref{Timo_weak_strong_Dirichlet_homogen}. Then,
 \begin{align}\label{lower_order_term2}
& \int\limits_0^\ell \left( \rho w_t(t,x)w(t,x)+I_\rho \psi_t(t,x) \psi(t, x)\right) dx\leq C_{DL} E(t),\;\;
\;  \mbox{ for case 1.-4. together with   \eqref{Timo_weak_strong_Dirichlet_homogen},  }  \;
\end{align}
and
 \begin{align}\label{lower_order_term21}
& \int\limits_0^\ell \left( \rho w_t(t,x)w(t,x)+I_\rho \psi_t(t,x) \psi(t, x)\right) dx\leq C_{NL} E_{\gamma,\delta}(t),
\;\;
\;  \mbox{for case 1.-4. together with   \eqref{Timo_weak_strong_Neumann_homogen},  }  \;
\end{align}
where
$$
C_{DL} = \max\{\sqrt{\rho},\sqrt{I_\rho}\} \sqrt{C_{D,K,EI}},
$$
$$
C_{NL} = \\max\{\sqrt{\rho},\sqrt{I_\rho}\} \sqrt{C_{N,K,EI}}.
$$

 \end{proposition}
 \begin{Proof}
 Indeed, when the boundary condition \eqref{Timo_weak_strong_Dirichlet_homogen} holds, we have
 \begin{align*}
 &\int\limits_0^\ell \left( \rho w_t(t,x)w(t,x)+I_\rho \psi_t(t,x) \psi(t, x)\right) dx
 \\ \leq \;&
 {1\over2}\sqrt{   C_{D,K,EI}}\max\{\sqrt{\rho},\sqrt{I_\rho}\} \int\limits_0^\ell  \left( \rho  w_t^2(t,x)  + I_\rho^2\psi_t^2(t,x) +  \frac{1}{C_{D,K,EI}}(w^2(t,x)+\psi^2(t,x))\right) dx.
 \end{align*}
 The analogous inequality holds for the case of boundary condition \eqref{Timo_weak_strong_Neumann_homogen}. Then  we can conclude the proof by Proposition \ref{Poincare_inequality} and \ref{Poincare_gd}.
 \end{Proof}
 We now come to the  main result of this section, the inverse inequality. We recall $\mu=\max(\mu_K,\mu_{EI})$.
 \begin{proposition}\label{Inverse}
\begin{enumerate}
\item[{\rm (i)}]
  Let $(w,\psi)$ be a classical solution of  \eqref{Timo_weak_strong_Neumann} with homogeneous Neumann conditions \eqref{Timo_weak_strong_Neumann_homogen} ($\gamma=\delta=0$) instead of \eqref{Timo_weak_strong_Neumann}$_5$. Then we have
  \begin{align}\label{inverseD}
 &\int\limits_0^T \ell \left(\rho w_t^2(t,\ell)+I_\rho \psi^2_t(t,\ell)\right)dt  \geq \left[ \left(2-\mu -2 {C}_{h}\right)T -4 C_{F}-  {\mu} C_{NL}\right]E(0).
 \end{align}
 \item[{\rm (ii)}] Let $(w,\psi)$ be a classical solution of \eqref{Timo_weak_strong_Neumann} with homogeneous Dirichlet conditions \eqref{Timo_weak_strong_Dirichlet_homogen} instead of \eqref{Timo_weak_strong_Neumann}$_5$. Then we have
      \begin{align}\label{inverseD1}
& \int\limits_0^T \ell \left(K(\ell)(w_x+\psi)^2(t,\ell)+ EI(\ell)\psi_x^2(t,\ell)\right)dt   \geq \left[ \left(2-\mu -2 {C}_{h}\right)T -4 C_{F}-  {\mu} C_{DL}\right]E(0).
 \end{align}

 \end{enumerate}
The constants $C_h,\; C_F,\; C_{NL},\; C_{DL}$
are defined in  \eqref{F_estimate_1}, \eqref{F_estimate_2}, \eqref{final_term},
\eqref{lower_order_term2}, and \eqref{lower_order_term21}, respectively.
 \end{proposition}
 \begin{Proof}
 We first treat the Neumann case. Recalling  \eqref{multiplier_identity_main} with $S=0$, we have
 \begin{align}\label{multiplier_identity_Dirichlet}
& \int\limits_0^T\ell  \left\{ \rho w_t^2(t,\ell)+I_\rho\psi_t^2(t,\ell) \right\} dt\notag\\
=\; & 2(F(T)-F(0))+\int\limits_0^T\int\limits_0^\ell \left\{ \rho w_t^2+I_\rho\psi_t^2 + (K-xK')(w_x+\psi)^2+(EI-xEI')\psi_x^2\right\}dx dt
\\& \quad\quad
+2\int\limits_0^T\int\limits_0^\ell x\left(-\rho  w_t \psi_t + K(w_x+\psi)\psi_x \right)dx d.\notag
\end{align}
We also recall the identity \eqref{multiplier_identity4} (with $S=0$) which we
multiply  by $\frac{\mu}{2}$ and add the result to \eqref{multiplier_identity_Dirichlet}. We obtain
\begin{align}\label{identity1}
& \int\limits_0^T\ell  \left\{ \rho w_t^2(t,\ell)+I_\rho\psi_t^2(t,\ell) \right\} dt\notag\\
=\; &(1-\frac{\mu}{2} )\int\limits_0^T\int\limits_0^\ell \left\{ \rho w_t^2+I_\rho \psi_t^2 + K(w_x+\psi)^2+EI\psi_x^2\right\}dx dt\\
& \quad\quad\quad\quad + \int\limits_0^T\int\limits_0^\ell \left\{ (\mu K-xK')(w_x+\psi)^2+(\mu EI-xEI')\psi_x^2\right\}dx dt+2(F(T)-F(0))\notag\\
&\quad\quad\quad\quad \quad\quad+\frac{\mu}{2}\int\limits_0^\ell \left( \rho w_t w +I_\rho \psi_t  \psi \right) dx|_0^T
+2\int\limits_0^T\int\limits_0^\ell x\left(-\rho  w_t \psi_t + K(w_x+\psi)\psi_x \right)dx dt.\notag
\end{align}
According to our assumption (H), the second integral is positive. Then,
\begin{align}\label{identity10}
& \int\limits_0^T\ell  \left\{ \rho w_t^2(t,\ell)+I_\rho\psi_t^2(t,\ell) \right\} dt
\ge\;
(2-\mu)\int\limits_0^T E(t) dt  +2(F(T)-F(0)) \\
&\quad\quad\quad\quad \quad\quad \quad \quad  +\frac{\mu}{2}\int\limits_0^\ell \left( \rho w_t w +I_\rho \psi_t  \psi \right) dx|_0^T
+2\int\limits_0^T\int\limits_0^\ell x\left(-\rho  w_t \psi_t + K(w_x+\psi)\psi_x \right)dx dt.\notag
\end{align}
 The remaining terms are estimated
by \eqref{final_term}, Proposition \ref{F_estimate} and \ref{lower_order_estimate}.

 Putting these estimates together, we arrive at \eqref{inverseD}. The Dirichlet case is treated in a similar way.
 \end{Proof}

 We finally treat the case of Robin conditions \eqref{Timo_weak_strong_Neumann_homogen}.
 This corresponds to changing  \eqref{identity1} to
 \begin{align}\label{identity11}
& \int\limits_0^T \ell\left\{ \rho w_t(t,\ell)^2+I_\rho \psi_t(t,\ell)^2+\frac{\gamma^2}{K(\ell)}w(t,\ell)^2+\frac{\delta^2}{EI(\ell)}\psi(t,\ell)^2\right\}dt.
\\
=\; &(1-\frac{\mu}{2} )\int\limits_0^T\int\limits_0^\ell \left\{ \rho w_t^2+I_\rho \psi_t^2 + K(w_x+\psi)^2+EI\psi_x^2\right\}dx dt\notag\\
&   + \int\limits_0^T\int\limits_0^\ell \left\{ (\mu K-xK')(w_x+\psi)^2+(\mu EI-xEI')\psi_x^2\right\}dx dt+2(F(T)-F(0))\notag\\
& +\frac{\mu}{2}\int\limits_0^\ell \left( \rho w_t w +I_\rho \psi_t  \psi \right) dx|_0^T
+2\int\limits_0^T\int\limits_0^\ell x\left(-\rho  w_t \psi_t + K(w_x+\psi)\psi_x \right)dx dt
+{\mu\over2}\int\limits_0^T  \left(\gamma   w^2(\ell)+ \delta \psi^2(\ell)\right)  dt.\notag
\end{align}
%
 As we now have to take the energy as \eqref{E-gd},
  we need to adjust for the extra  boundary terms in \eqref{inverseD} to obtain
 \begin{align}\label{Inverse_Robin_0}
&\int\limits_0^T\ell  \Big[ \rho w_t^2(t,\ell)+I_\rho \psi_t^2 (t,\ell)
+\Big(\frac{\gamma}{K(\ell)}+2-\mu-2{C}_{h}\Big)\gamma w^2(t,\ell)
 +\Big(\frac{\delta}{EI(\ell)}+2-\mu-2{C}_{h}\Big)\delta\psi^2(t,\ell)\Big]dt
 \\
 \geq \;\;  & \big(2-\mu-2{C}_{h}\big)\int\limits_0^T E_{\gamma,\delta}(t)dt -
(4 C_{F}+  {\mu} C_{NL})E(0)+{\mu\over2}\int\limits_0^T  \left(\gamma   w^2(\ell)+ \delta \psi^2(\ell) \right)  dt.
\notag
 \end{align}
  It obvious that the direct inequalities are not effected by the change of the energy expression. Therefore, we conclude:

 \begin{proposition}\label{Inverse_Robin}
 Let $(w,\psi)$ be smooth solutions of \eqref{Timo_weak_strong_Neumann} with Robin conditions \eqref{Timo_weak_strong_Neumann_homogen} ($\gamma,\delta\geq 0$) instead of \eqref{Timo_weak_strong_Neumann}$_5$.  
 Then, we have the following inverse inequality
 \begin{align} \label{eq-eta}
 &\int\limits_0^T \ell
 \left[ \rho w_t^2(t,\ell)+I_\rho \psi_t^2(t,\ell)
 +\eta_1\gamma w^2(t,\ell)+\eta_2\delta\psi^2(t,\ell)\right]dt
 \\
 \geq\;\;&
\left[ \left(2-\mu -2 {C}_{h}\right)T -4 C_{F}-  {\mu} C_{NL}\right]E_{\gamma,\delta}(0).\notag
 \end{align}
 where
 \[
\eta_1:= \frac{\gamma}{K(\ell)}+2-\mu-2{C}_{h}, \; \eta_2:=\frac{\delta}{EI(\ell)}+2-\mu-2{C}_{h}.  \]
 \end{proposition}

 We can use the compactness-uniqueness argument in \cite{Komornik1989} Proposition 6.18 to conclude
 \begin{proposition}\label{Invers_Robin_improved}
 We assume the constant defined in \eqref{final_term} satisfies
\begin{equation}\label{small}
\mu+ 2{C}_{h}<2.
\end{equation}
 Under the conditions of Proposition \ref{Inverse_Robin},
  for $T>0$ sufficiently large, there exist a constant $C>0$ such that
 \begin{align}\label{Inverse_Robin_1}
\int\limits_0^T\ell
 \left\{ \rho w_t(t,\ell)^2+I_\rho \psi_t(t,\ell)^2\right\}dt\geq C E_{\gamma,\delta}(0),\;\; \forall \Big((w^0,\psi^0),(w^,\psi^1)\Big)\in V^2_{K,EI}\times V^1_{K,EI}.
 \end{align}
 Moreover, in the weakly degenerate case, we can  provide an explicit constant $C_w>0$ such that \eqref{Inverse_Robin_1} holds:
\begin{align}\label{constant_robin}
 C_w:=& \frac{1}{\left(1+\frac{2}{\ell}\max\{\frac{\eta_1\gamma}{\rho},
 \frac{\eta_2\delta}{I_\rho}\}\right)}
 \left[ \left(2-\mu -2{C}_h\right) T
 -4  C_{F}- {\mu}C_{NL}
 -4\max\Big\{\eta_1\gamma \Big\|{1\over K}\Big\|_{L^1},\eta_2\delta \Big\|{1\over EI}\Big\|_{L^1}\Big\}\right].
 \end{align} 
 \end{proposition}
 \begin{Proof}
    By  a standard one by contradiction (see \cite[Prop. 6.18]{Komornik1989}), we can get that there exists positive constant $\tilde{C}$ such that
     $$\int\limits_0^\ell ( w^2(t,\ell)+ \psi^2(t,\ell)) dt
     \le  \tilde{C}\int\limits_0^\ell(w_t^2(t,\ell)+ \psi_t^2(t,\ell))dt.
     $$
      Therefore, no further information on the constant $C$ is available. Thus, the observability time is not given explicitly in this case.

      In the case of weak degeneracy, we can estimate as follows:
     \[
     |w(0,\ell)|^2
     =\Big(\int\limits_0^\ell w_x(0,s)ds\Big)^2\leq \Big\|{1\over K}\Big\|_{L^1}|w(0,\cdot)|^2_{1,K}\leq 2 \Big\|{1\over K}\Big\|_{L^1}E_{\gamma,\delta}(0).
     \]
     Hence,
     \[
     |w(t,\ell)|^2=\Big|\int\limits_0^t w_s(s,\ell)ds+w(0,\ell)\Big|^2\leq 2\int\limits_0^T w_t^2(t,\ell)dt +2 w^2(0,\ell)\leq 2\int\limits_0^T w_t^2(t,\ell)dt + 4 \Big\|{1\over K}\Big\|_{L^1} E_{\gamma,\delta}(0).     \]
     The analogous statement holds for $\psi$. Thus,
     \begin{align*}
       & \int\limits_0^T \eta_1 \gamma w^2(t,\ell)+\eta_2\delta \psi^2(t,\ell) dt\\
   \leq \;   & \frac{2}{\ell}\max\Big\{\frac{\eta_1\gamma}{\rho},\frac{\eta_2\delta}{I_\rho}\Big\}\int\limits_0^T \ell \left( \rho w_t^2(t,\ell)+I_\rho \psi_t^2(t,\ell)\right)dt
       +4\max\Big\{\eta_1\gamma\Big\|{1\over K}\Big\|_{L^1},\eta_2\delta \Big\|{1\over EI}\|_{L^1}\Big\}E_{\gamma,\delta}(0).
    \end{align*}
   By Proposition \ref{Inverse_Robin}, this implies
    \begin{align*}
 &\left(1+\frac{2}{\ell}\max\Big\{\frac{\eta_1\gamma}{\rho},\frac{\eta_2\delta}{I_\rho}\Big\}\right)\int\limits_0^T \ell \left( \rho w_t(t,\ell)^21+I_\rho \psi_t(t,\ell)^2\right)dt
 \\
 \geq\; & \left[ \left(2-\mu -2{C}_h\right) T
 -4  C_{F}- {\mu}C_{NL}
 -4\max\Big\{\eta_1\gamma \Big\|{1\over K}\Big\|_{L^1},\eta_2\delta \Big\|{1\over EI}\Big\|_{L^1}\Big\}\right]E_{\gamma,\delta}(0).
 \end{align*}
 \end{Proof}
 We now argue as in \cite{Komornik1989} and prove
 \begin{proposition}\label{Invers_Robin_final}
  Let the shifted energy norm be given by
  $$E^{-\frac{1}{2}}((w^0,\psi^0),(w^1,\psi^1):=
  \frac{1}{2}\{\|(w^0,\psi^0\|_H^2+\||(w^1,\psi^1)\||^2_{(V^1_{K,EI})^*}\}.$$ Then
 \begin{align}\label{Inverse_Robin_3}
\int\limits_0^T\ell
 \left\{ \rho w^2(t,\ell)+I_\rho \psi^2(t,\ell)\right\}dt\geq C E^{-\frac{1}{2}}(0),\;\; \forall \Big((w^0,\psi^0),(w^1,\psi^1)\Big)\in V^1_{K,EI}\times H,
 \end{align}
  with the same $C$ as in \eqref{constant_robin}.
 \end{proposition}

 \section{Observability}\label{observability}
 Exact boundary observability in time $T>T_0$ of the Timoshenko systems \eqref{Timo_weak_strong_Neumann} with boundary conditions \eqref{Timo_weak_strong_Dirichlet_homogen} or \eqref{Timo_weak_strong_Neumann_homogen} relates to the question as to whether the initial states  can be reconstructed by observed data. In our situation, the homogeneous boundary conditions are \eqref{Timo_weak_strong_Dirichlet_homogen} and \eqref{Timo_weak_strong_Neumann_homogen} while the corresponding observed data are $K(\ell)(w_x(t,\ell)+\psi(t,\ell), EI(\ell)\psi_x(t,\ell), t\in (0,T)$ and $w_t(t,\ell), \psi_t(t,\ell), t\in (0,T)$ (or $w(t,\ell), \psi(t,\ell), t\in (0,T)$), respectively. The lower bound $T_0$ of observability time depends on the speed of propagation of the two waves in terms of $w$ and $\psi$. It is known, see e.g. \cite{Taylor00}, that the two velocities are given by
 \[
 v_1(x):=\sqrt{\frac{K(x)}{\rho}}, \;\; v_2(x):=\sqrt{\frac{EI(x)}{I_\rho}}. \]
 Moreover, the time a signal needs to travel cross the domain is
 \[
T_1:=\int\limits_0^\ell \frac{1}{v_1(x)}dx=\int\limits_0^\ell \sqrt{\frac{\rho}{K(x)}}dx,\;\; T_2:= \int\limits_0^\ell \sqrt{\frac{I_\rho}{EI(x)}}dx. \]
\begin{remark}\label{observability_time}
   We remark that in the case of the exemplary degeneracies $K(x)=Kx^{\mu_K}$, $EI(x)=EI x^{\mu_{EI}}$, the travel times can be computed explicitly.
   \[
   T_{1,\mu_K}=\sqrt{\frac{\rho}{K}}\frac{2}{2-\mu_K}\ell^{\frac{2-\mu_K}{2}},\;\; T_{2,\mu_{EI}}=\sqrt{\frac{I_\rho}{EI}}\frac{2}{2-\mu_{EI}}\ell^{\frac{2-\mu_{EI}}{2}}.
   \]
   In particular, for $\mu_K= \mu_{EI} =0$, non-degenerate case, we obtain
   \[
   T_{1,0}=\sqrt{\frac{\rho}{K}}\ell, \;\; T_{2,0}=\sqrt{\frac{I_\rho}{EI}}\ell.
   \]
 In the case of one-sided observability (exact controllability) for \it{non-degenerate} Timoshenko beams, we have
 \[ T_0=2\max(T_{1,0}, T_{2,0}).\]
 In the case of degeneration, the time $T_0$ is provided by the inverse inequalities below.
\end{remark}
 The observability requirement is then equivalent to the inequalities provided in Proposition \ref{Inverse}, \ref{Inverse_Robin} and \ref{Invers_Robin_improved} as they provide the left invertibility of the observation map. Therefore, Proposition \ref{Inverse}, \ref{Inverse_Robin} and \ref{Invers_Robin_improved} provides the argument. The only point to add is that the inequalities there have been established for smooth solutions. However, given the semi-group approach (or equivalently the approach using cosine and sine propagators) as well as the approach by transposition for very weak solutions in the Dirichlet case, can be employed using the fact that smooth solutions are dense in the corresponding mild or weak solutions. Thus, the corresponding inequalities can be extended to such solutions. We can, therefore, conclude:


\begin{theorem}\label{observability_cases}
  Assume that the smallness assumption \eqref{small} holds.  We distinguish three cases:
\begin{description}
\item[{\rm(i)} {\it Dirichlet case:}]
    System \eqref{Timo_weak_strong_Neumann} with \eqref{Timo_weak_strong_Dirichlet_homogen} instead of \eqref{Timo_weak_strong_Neumann}$_5$ is observable within finite energy solutions $((w^0,\psi^0),\-(w^1,\psi^1)) \in \mathcal{V}_{D,w/s}\times {\mathcal H}$ in time $T$, with
    \begin{equation}\label{Dirichlet_time}
        T > \left(2-\mu -2 {C}_{h}\right)^{-1} \left(4 C_{F}-  {\mu} C_{DL}\right).
    \end{equation}

   \item[{\rm(ii)} {\it Robin case:}] System \eqref{Timo_weak_strong_Neumann} with \eqref{Timo_weak_strong_Neumann_homogen}instead of \eqref{Timo_weak_strong_Neumann}$_5$ and $\gamma,\delta > 0$, in case of weak degeneracy,  is observable within shifted energy solutions $((w^0,\psi^0),\-(w^1,\psi^1))\in  {\mathcal H}\times (\mathcal{V}_{R,w})^*$ in time $T$, with
    \begin{align}\label{constant_robin_time}
    T>\left(2-\mu -2{C}_h\right)^{-1} \left(
 4  C_{F}+{\mu}C_{NL}
 +4\max\Big\{\eta_1\gamma \Big\|{1\over K}\Big\|_{L^1},\eta_2\delta \Big\|{1\over EI}\Big\|_{L^1}\Big\}\right).
    \end{align}

    In the strong degenerate case, the system is observable within shifted energy solutions $((w^0,\psi^0),\-(w^1,\psi^1))\in  {\mathcal H}\times (\mathcal{V}_{R,s})^*$ in time $T$.
     In this case, however, we can only say that $T>0$ has to be sufficiently large.


\item[{\rm(iii)} {\it Neumann case:}] System  \eqref{Timo_weak_strong_Neumann} with \eqref{Timo_weak_strong_Neumann_homogen}instead of \eqref{Timo_weak_strong_Neumann}$_5$ and    $\gamma=\delta=0$,
    in the weak degenerate case, is observable within finite energy solutions $((w^0,\psi^0),\-(w^1,\psi^1)) \in \mathcal{V}_{N,w/s}\times {\mathcal H}$ in time $T$, with
    \begin{equation} \label{time_neu}
        T > \left(2-\mu -2 {C}_{h}\right)^{-1} \left(4 C_{F}-  {\mu} C_{NL}\right).
    \end{equation}
   In the strong degenerate case, the system is observable within shifted energy solutions, where the data are orthogonal to constants, and the time $T>0$ satisfies \eqref{time_neu}.
\end{description}
The constants $C_h,\; C_F,\;   C_{DL},\;   C_{NL},\; \eta_1,\; \eta_2,\; C_w$
 are defined in  \eqref{F_estimate_1}, \eqref{F_estimate_2}, \eqref{final_term},
\eqref{lower_order_term2}, \eqref{lower_order_term21}, \eqref{eq-eta} and \eqref{constant_robin}, respectively.
    \end{theorem}

\begin{remark}
     We notice that the time $T>T_0$, the lower bound for observability in the non-degenerated case given by Remark \ref{observability_time} can be revealed in case $\mu_K=\mu_{EI}=0$  for very thin beams. We remark that using the method of characteristics as in \cite{LagneseLeugeringSchmidt1993b} or using extensions as in \cite{Taylor00}, no constraint on the thickness is present. In case of degeneration, we have this time as lower bound, too. See the estimates \eqref{F_estimate_1}, \eqref{F_estimate_2} in each particular case. Clearly, as $\mu\rightarrow 2$, $T_0\rightarrow \infty$. This suggests that we loose observability in this case. The explicit relation of the constants with respect to the geometrical and physical parameters allows one to determine the precise observability times.
 \end{remark}

 \section{Controllability}\label{controllability}
   \subsection{Solutions in the sense of transposition}
 \subsubsection{Dirichlet case}
 For the degenerate control system,  we develop the meaning of transposition.
 First, for the system \eqref{Timo_weak_strong_Neumann}$_{1,2,3,4}$ with \eqref{Timo_weak_strong_Dirichlet}, we assume that $  f_1,\;f_2 \in L^2_{\text{loc}}(0, \infty).$
  Take the pair $(v,\phi)$ as a solution of the backward system 
\begin{align}\label{Timo_weak_strong_Dirichlet_homogen_final}
& \rho v_{tt}-\big(K(x)(v_x+\phi)\big)_x=0,\; (t,x)\in Q,\notag\\
& I_\rho \phi_{tt}-(EI(x)\phi_x)_x+K(x)(v_x+\phi)=0,\; (t,x)\in Q,\notag\\
& v(t,0)=0, \; \text{ if } \mu_K\in[0,1),\; \lim\limits_{x\rightarrow 0^+}K(x)(v_x+\phi)(t,x)=0\;  \text{ if } \; \mu_K\in [1,2),\; t\in (0,T),\\
& \phi(t,0)=0, \; \text{ if } \mu_{EI}\in[0,1),\; \lim\limits_{x\rightarrow 0^+}EI(x)\phi_x(t,x)=0\;  \text{ if } \; \mu_EI\in [1,2),\; t\in (0,T),\notag\\
& v(t,\ell)=0, \; \phi(t,\ell)=0, \; t\in (0,T),\notag\\
& v(T,x)=v^0_T(x), \; v_t(T,x)=v^1_T(x), \; \phi(T,x)=\phi^0_T(x), \; \phi_t(T,x)=\phi^1_T(x), \; x\in (0,\ell).\notag
\end{align}
Indeed, according to the well-posedness result for   homogeneous systems, we have the existence and uniqueness of a pair $(v,\phi)\in C^1([0, \infty); {\cal H}) \cap C([0, \infty); \mathcal{V}_{D,w/s})$ satisfying \eqref{Timo_weak_strong_Dirichlet_homogen_final} with initial data
\begin{equation}\label{final_data}
(v_T^0,\phi_T^0)\in \mathcal{V}_{D,w/s}(0,\ell), \;\;  (v_T^1,\phi_T^1)\in {\mathcal H}. 
\end{equation}
We first take smooth data and multiply the equations \eqref{Timo_weak_strong_Neumann}$_{1,2}$ by the corresponding $v, \phi$, respectively, and use \eqref{Timo_weak_strong_Neumann}$_{3,4}$, \eqref{Timo_weak_strong_Dirichlet}. Then we obtain after integration by parts
\begin{align}\label{transposition_Dirichlet1}
0=&\int\limits_0^T\int\limits_0^\ell\left\{\left( \rho w_{tt}-(K(w_x+\psi))_x \right)v+\left( I_\rho\psi_{tt}-(EI\psi_x)_x +K(w_x+\psi)\right) \phi \right\}dx dt\\
& =\int\limits_0^\ell \left( \rho w_t v-\rho w v_t +I_\rho \psi_t \phi -I_\rho \psi \phi_t\right) dx|_0^T+\int\limits_0^T\left(f_2 EI(\ell)\phi_x( \ell)+f_1 K(\ell)(v_x+\phi)( \ell)\right)dt.\notag
\end{align}
This can be written as
\begin{align}\label{eq28}
&\int\limits_0^T\left(f_2(t)EI(\ell)\phi_x(t,\ell)+f_1(t)K(\ell)(v_x+\phi)(t,\ell)\right)dt\\
   =\;  & \Big(\begin{pmatrix}
       w_t(T)\\ \psi_t(T)
   \end{pmatrix},  \begin{pmatrix}
       v_T^0\\ \phi_T^0\end{pmatrix}\Big)_{{\mathcal H}\times {\mathcal H}}-\Big(\begin{pmatrix}
       w(T)\\ \psi(T)
   \end{pmatrix},  \begin{pmatrix}
       v_T^1\\ \phi_T^1 \end{pmatrix}\Big)_{{\mathcal H}\times {\mathcal H}}\notag\\
       & \qquad\qquad\qquad - \Big(\begin{pmatrix}
       w_t(0)\\ \psi_t(0)
   \end{pmatrix},  \begin{pmatrix}
       v(0)\\ \phi(0)\end{pmatrix}\Big)_{{\mathcal H}\times {\mathcal H}}+\Big(\begin{pmatrix}
       w(0)\\ \psi(0)
   \end{pmatrix},  \begin{pmatrix}
       v_t(0)\\ \phi_t(0)   \end{pmatrix}\Big)_{{\mathcal H}\times {\mathcal H}}.\notag\end{align}

 Now we define \( (\mathcal{V}_{D,w/s})^*\) as the dual space of \( \mathcal{V}_{D,w/s}  \) with respect to the pivot space \( {\mathcal H} \). Then, thanks to Proposition \ref{Poincare_inequality},  the operator  \( A  \)
 with domain $ D(A)_{D,w/s}$ is an isomorphism from \(\mathcal{V}_{D,w/s} \) onto \( (\mathcal{V}_{D,w/s})^* \).
 Then it follows from \eqref{eq28} that
\begin{align}\label{transposition_Dirichlet2}
&\int\limits_0^T\left(f_2(t)EI(\ell)\phi_x(t,\ell)+f_1(t)K(\ell)(v_x+\phi)(t,\ell)\right)dt\\
&\qquad + \Big\langle\begin{pmatrix}
       w^1\\ \psi^1
   \end{pmatrix},\begin{pmatrix}
       v(0)\\ \phi(0)\end{pmatrix}\Big\rangle_{(\mathcal{V}_{D,w/s})^*\times \mathcal{V}_{D,w/s}}-\Big\langle\begin{pmatrix}
       w^0\\ \psi^0
   \end{pmatrix},\begin{pmatrix}
       v_t(0)\\ \phi_t(0)   \end{pmatrix}\Big\rangle_{\mathcal{H}\times \mathcal{H}} \notag\\
       = \; & \Big\langle\begin{pmatrix}
       w_t(T)\\ \psi_t(T)
   \end{pmatrix},\begin{pmatrix}
       v_T^0\\ \phi_T^0\end{pmatrix}\Big\rangle_{(\mathcal{V}_{D,w/s})^*\times \mathcal{V}_{D,w/s}}-\Big\langle\begin{pmatrix}
       w(T)\\ \psi(T)
   \end{pmatrix},\begin{pmatrix}
       v_T^1\\ \phi_T^1 \end{pmatrix}\Big\rangle_{\mathcal{H}\times \mathcal{H}}\notag\\
       & \qquad\qquad\qquad  \notag
\end{align}
for all $T>0$ and final data \eqref{final_data}.

It is known by Proposition \ref{conservation} that the energy for system \eqref{Timo_weak_strong_Dirichlet_homogen_final} is conserved. Moreover, from Proposition \ref{upper_Dirichlet}, we have
$$
\int\limits_0^T \ell \left(K(\ell)v_x^2(t,\ell)+ EI(\ell) \phi^2_x(t,\ell)\right)dt\leq C_{BC}E(0) =C_{BC}E(T).
$$
Therefore, the left hand side defines a continuous linear form on $(v_T^0,\phi_T^0, v_T^1,\phi_T^1)\in  \mathcal{V}_{D,w/s}\times\mathcal{H}$ and so is the right hand side.
  Therefore, there is a unique solution of \eqref{Timo_weak_strong_Neumann}$_{1,2,3,4,6}$ with \eqref{Timo_weak_strong_Dirichlet}  by transposition
$$
(w,\psi) \in C^1([0, \infty); (\mathcal{V}_{D,w/s})^*)
\cap C([0, \infty);{\cal H} ).
$$
%
 \subsubsection{Robin case}
 In this case $(\gamma,\delta>0)$, standard semi-group or  cosine propagator theory gives us the regularity result
 $(w,\psi)\in C^1(0,T; {\mathcal H})\cap C(0,T;\mathcal{V}_{R,w/s})$.
With respect to the notion of solution in the sense of transposition, in this case, we consider system \eqref{Timo_weak_strong_Dirichlet_homogen_final} but now with boundary conditions
\begin{align}\label{Timo_weak_strong_Neuman_homogen_final}
K(\ell)(v_x(t,\ell)+\phi(t,\ell))+\gamma v(t,\ell)=0, \; EI(\ell)\phi_x(t,\ell)+\delta \phi(t,\ell)=0, \; t\in (0,T),
\end{align}
with $\gamma, \delta> 0$.
For $(w,\psi)$ satisfies \eqref{Timo_weak_strong_Neumann}$_{1,2,3,4,5}$, a similar calculation as in the Dirichlet case leads to:
 \begin{align}\label{transposition_Neumann}
   & \int\limits_0^T \left(f_1(t)v(t,\ell)+f_2(t) \phi(t,\ell) \right)dt \\
   & =\Big\langle\begin{pmatrix}
       w_t(T)\\ \psi_t(T)
   \end{pmatrix},\begin{pmatrix}
       v_T^0\\ \phi_T^0\end{pmatrix}\Big\rangle_{\mathcal{H}\times \mathcal{H}}-\Big\langle\begin{pmatrix}
       w(T)\\ \psi(T)
   \end{pmatrix},  \begin{pmatrix}
       v_T^1\\ \phi_T^1 \end{pmatrix}\Big\rangle_{\mathcal{V}_{R,w/s}\times (\mathcal{V}_{R,w/s})^*}\notag\\
       & \qquad\qquad\qquad - \Big\langle\begin{pmatrix}
       w^1\\ \psi^1
   \end{pmatrix},\begin{pmatrix}
       v(0)\\ \phi(0)\end{pmatrix}\Big\rangle_{{\mathcal H}\times {\mathcal H}}+\Big\langle\begin{pmatrix}
       w^0\\ \psi^0
   \end{pmatrix}, \begin{pmatrix}
       v_T(0)\\ \phi_t(0)   \end{pmatrix}\Big\rangle_{\mathcal{V}_{R,w/s)}\times (\mathcal{V}_{R,w/s})^*},\notag \end{align}
      for smooth data. We obtain again the existence of a unique solution in the sense of transposition.
      \subsubsection{The Neumann case}
      In this case $\gamma=\delta=0$. For weak degeneracy, everything stays as in the Robin-case. Only in the strongly degenerate case, we have to replace $\mathcal{V}_{R,s}$ by $\mathcal{V}_{N,s}$, and hence we work in the quotient space. The same argument as above shows that there exists a unique solution in the sense of transposition.

\subsection{Exact controllability}
 We will show exact controllability using the arguments in the HUM method. The idea is to use the direct and indirect inequalities established in Section \ref{inequalities} and Section \ref{observability} in order to show  that given final data $X^T=(w_T^0,\psi_T^0,w_T^1, \psi_T^1)$ and $\tilde{X}^T=(\tilde{w}_T^0,\tilde{\psi}_T^0,\tilde{w}_T^1, \psi_T^1)$ for the homogeneous backwards running (adjoint) systems \eqref{Timo_weak_strong_Dirichlet_homogen_final} or \eqref{Timo_weak_strong_Neuman_homogen_final}, the bilinear forms
 \begin{align}\label{HUM-forms}
     \Lambda_D(X^T,\tilde{X}^T)&:= \int\limits_0^T \big( K(\ell)(w_x+\psi)(t,\ell)(\tilde{w}_x+\tilde{\psi})(t,\ell) +EI(\ell)\psi_x(t,\ell)\tilde{\psi}_x(t,\ell)\big)dt, \;\; \forall X^T,\tilde{X}^T \in \mathcal{V}_{D,w/s}\times \mathcal{H},\\
 \Lambda_R(X^T,\tilde{X}^T)&:= \int\limits_0^T (w(t,\ell)\tilde{w}(t,\ell) + \psi(t,\ell)\tilde{\psi}(t,\ell) )dt, \;\; \forall X^T,\tilde{X}^T \in  \mathcal{H}\times (\mathcal{V}_{D,w/s})^*,
 \end{align}
 are continuous and coercive, respectively. We note that due to the time reversibility of the Timoshenko beam system, the problem of exact controllability is equivalent to that of null controllability. This means it is sufficent to consider the situation that $(w(T),\psi(T)), (w_t(T),\psi_t(T))=0$.  We look into the cases separately.

 \subsubsection{Dirichlet case}
In this case according Proposition \ref{Inverse}(ii) or Theorem \ref{observability_cases} (i), we have that $\Lambda_D$ is continuous and coercive for $T$ satisfying \eqref{Dirichlet_time}. We define the continuous linear functional
\begin{equation}
   \mathcal{L}_D(X^T):= \Big\langle\begin{pmatrix}
       w^1\\ \psi^1
   \end{pmatrix},\begin{pmatrix}
       v(0)\\ \phi(0)\end{pmatrix}\Big\rangle_{(\mathcal{V}_{D,w/s)})^*\times \mathcal{V}_{D,w/s}}-\Big\langle\begin{pmatrix}
       w^0\\ \psi^0
   \end{pmatrix},\begin{pmatrix}
       v_t(0)\\ \phi_t(0)   \end{pmatrix}\Big\rangle_{\mathcal{H}\times \mathcal{H}},\;\; \forall X^T\in \mathcal{V}_{D,w/s}\times \mathcal{H}.\end{equation}
 As $\Lambda_D$ is a continuous and coercive bilinear form on 
 $(\mathcal{V}_{D,w/s}\times \mathcal{H})^2$, the variational equation
 \begin{equation}\label{HUM_D}
      \Lambda_D(X^T,\tilde{X}^T)=-\mathcal{L}_D(\tilde{X}^T), \; \; \forall \tilde{X}^T\in \mathcal{V}_{D,w/s}\times \mathcal{H}
 \end{equation}
 has a unique solution $X^T\in \mathcal{V}_{D,w/s}\times \mathcal{H}$.

 We now set $f_1(t):= w_x(t,\ell)+\psi(t,\ell)$ and $f_2(t,\ell)=\psi_x(t,\ell)$. Then \eqref{HUM_D} reads
 \begin{align}
     &\int\limits_0^T f_1(t)K(\ell)(\tilde{w}_x(t,\ell)+\tilde{\psi}(t,\ell))+f_2(t)EI(\ell)\tilde{\psi}(t,\ell)dt=\Lambda_D(X^T, \tilde{X}^T)\notag\\
     & \qquad = -\Big\langle\begin{pmatrix}
       w^1\\ \psi^1
   \end{pmatrix},\begin{pmatrix}
       v(0)\\ \phi(0)\end{pmatrix}\Big\rangle_{(\mathcal{V}_{D,w/s)})^*\times \mathcal{V}_{D,w/s}}+\Big\langle\begin{pmatrix}
       w^0\\ \psi^0
   \end{pmatrix},\begin{pmatrix}
       v_t(0)\\ \phi_t(0)   \end{pmatrix}\Big\rangle_{\mathcal{H}\times \mathcal{H}}, \; \forall \tilde{X}^T\in\mathcal{V}_{D,w/s}\times \mathcal{H} .\end{align}
       On the other side, $f_1,f_2$ generate a solution $(w,\psi)$  in the sense of transposition, hence, by \eqref{transposition_Dirichlet2},
   \begin{align}
        &\int\limits_0^T f_1(t)K(\ell)(\tilde{w}_x(t,\ell)+\tilde{\psi}(t,\ell))+f_2(t)EI(\ell)\tilde{\psi}(t,\ell)dt\notag \\
     & \qquad = -\Big\langle\begin{pmatrix}
       w^1\\ \psi^1
   \end{pmatrix},\begin{pmatrix}
       v(0)\\ \phi(0)\end{pmatrix}\Big\rangle_{(\mathcal{V}_{D,w/s)})^*\times \mathcal{V}_{D,w/s}}+\Big\langle\begin{pmatrix}
       w^0\\ \psi^0
   \end{pmatrix},\begin{pmatrix}
       v_t(0)\\ \phi_t(0)   \end{pmatrix}\Big\rangle_{\mathcal{H}\times \mathcal{H}}\\
       & \qquad\qquad +\Big\langle\begin{pmatrix}
       w_t(T)\\ \psi_t(T)
   \end{pmatrix},\begin{pmatrix}
       v_T^0\\ \phi_T^0\end{pmatrix}\Big\rangle_{(\mathcal{V}_{D,w/s})^*\times \mathcal{V}_{D,w/s}}-\Big\langle\begin{pmatrix}
       w(T)\\ \psi(T)
   \end{pmatrix},\begin{pmatrix}
       v_T^1\\ \phi_T^1 \end{pmatrix}\Big\rangle_{\mathcal{H}\times \mathcal{H}}, \; \forall \tilde{X}^T\in \mathcal{V}_{D,w/s}\times \mathcal{H}. \notag  \end{align}
       Then, by the uniqueness of the solutions $(w,\psi)$ in the sense of transposition, we obtain $w(T)=0, w_t(T)=0, \psi(T)=0, \psi_t(T)=0$, as required.
       We thus have proved
       \begin{theorem}\label{controllability_D}
           Assume \eqref{small} holds and $T$ satisfies \eqref{Dirichlet_time}p. Then all initial and final states in $\mathcal{H}\times(\mathcal{V}_{D,w/s})^*$, the system   \eqref{Timo_weak_strong_Neumann} with \eqref{Timo_weak_strong_Dirichlet} instead of \eqref{Timo_weak_strong_Neumann}$_5$  is exactly controllable by controls $f_1,f_2\in L^2(0,T).$
       \end{theorem}

 \subsubsection{Robin case}
 We do not repeat the arguments here, but rather recognize that according to Theorem \ref{observability_cases} ii.) the form $\Lambda_R$ is continuous and coercive in $(\mathcal{H}\times \mathcal{V}_{R,w/s})^*$ for $T$ satisfying \eqref{constant_robin_time}. The corresponding solutions $(w,\psi)$ in the sense of transposition are given in the original energy space and, therefore, we obtain.
  \begin{theorem}\label{controllability_R}
           Assume $T$ satisfies \eqref{constant_robin_time} and \eqref{small} holds. Then all initial and final states in $\mathcal{V}_{D,w/s}\times\mathcal{H}$, the system   \eqref{Timo_weak_strong_Neumann} is exactly controllable by controls $f_1,f_2\in L^2(0,T).$
       \end{theorem}

       \subsubsection{Neumann case}
       The case of Neumann controls can be treated in the same way as the Robin case as long as weak degeneracy occurs. In the case of strong degeneracy, we work with the corresponding quotient space. The details are left to the reader.
 \section{Stabilzability}\label{stabilizability}
We consider the degenerated Timoshenko system  \eqref{Timo_weak_strong_Neumann} with state and velocity boundary feedback controls $f_1(t)=-\alpha w_t(t,\ell), f_2(t)=-\beta \psi_t(t,\ell)$, i.e.
\begin{align}\label{Timo_feedback}
& K(\ell)(w_x(t,\ell)+\psi(t,\ell))+\gamma w(t,\ell)=-\alpha w_t(t,\ell), \; EI(\ell)\psi(t,\ell)+\delta \psi(t,\ell)=-\beta \psi(t,\ell), \; t\in (0,T).
\end{align}
The corresponding modified energy is defined by \eqref{E-gd},
where $\gamma,\delta,\alpha,\beta\geq 0$ are given constants. In particular, we choose $\alpha,\beta >0$. The feedback law  represents a general state and velocity boundary feedback law. By simple calculations, we can show:
\begin{proposition}\label{dissipation}
Let $(w,\psi)$ be the solution of \eqref{Timo_weak_strong_Neumann} with \eqref{Timo_feedback}. Then
\begin{equation}\label{dissipative}
\frac{d}{dt}E_{\gamma,\delta}(t)= -\alpha w_t(t,\ell)^2-\beta \psi_t(t,\ell)^2, \; t\geq 0.
\end{equation}
Thus, the system \eqref{Timo_weak_strong_Neumann} with \eqref{Timo_feedback} instead of \eqref{Timo_weak_strong_Neumann}$_5$ is dissipative.
\end{proposition}
Well-posedness of the feedback control system in the finite energy space follows along the ideas in Subsection \ref{well-posendess} or the lines of \cite{AlabauCannarsaLeugering2017} for strings and, due to space limitations, the details are left to the reader. The question is as to whether we can show exponential decay of the energy. To this end, we reconsider the multiplier identity \eqref{multiplier_identity_main}. We recall the estimate \eqref{final_term}.
 Moreover, as in Proposition \ref{F_estimate}, we now have
\[
 F(T)-F(S)\leq  C_{F}(E_{\gamma,\delta}(T)+E_{\gamma,\delta}(S)) \]
 and as in Proposition \ref{lower_order_estimate}, we obtain
  \begin{align*}
& \int\limits_0^\ell \left( \rho w_t(t,x)w(t,x)+I_\rho \psi_t(t,x) \psi(t, x)\right)|_S^T dx\leq C_{NL}(E_{\gamma,\delta}(T)+E_{\gamma,\delta}(S)).
\end{align*}
We multiply \eqref{multiplier_identity3} by $\frac{\mu}{2}$ and add the result to \eqref{multiplier_identity_main} to obtain
\begin{align}\label{identity1_ST}
& \int\limits_S^T\ell  \left\{ \Big(\rho+\frac{\alpha^2}{K(\ell)}\Big)w_t^2(\ell)+\Big(I_\rho+\frac{\beta^2}{EI(\ell)}\Big)\psi_t^2(\ell) \right\} dt +\int\limits_S^T\left(\left(\frac{\ell}{K(\ell)}\gamma-\frac{\mu}{2}\right) \gamma w^2(\ell)\right.
\\
&\left. \quad 
+\left(\frac{\ell}{EI(\ell)}\delta-\frac{\mu}{2}\right)\delta\psi^2(\ell)\right) dt+\int\limits_S^T\left(\alpha\left(\frac{2\ell\gamma}{K(\ell)}-\frac{\mu}{2}\right) w_t(\ell)w(\ell)+\delta\left(\frac{2\ell\delta}{EI(\ell)}-\frac{\mu}{2}\right)\psi_t(\ell) \psi(\ell)\right) dt \notag
\end{align}
\begin{align*}
=\; & \frac{1}{2}\int\limits_S^T\int\limits_0^\ell\left\{(2-\mu) \left( \rho w_t^2+I_\rho \psi_t^2\right) + \left(2(K-xK')+\mu K\right))(w_x+\psi)^2+\left(2(EI-xEI'+\mu EI)\right)\psi_x^2\right\}dx dt\notag\\
&  \quad  +2(F(T)-F(S))+\frac{\mu}{2}\int\limits_0^\ell \left( \rho w_t w+I_\rho \psi_t \psi\right) dx|_S^T+2\int\limits_S^T\int\limits_0^\ell x\left(\rho  w_t \psi_t + K(w_x+\psi)\psi_x \right)dx dt\notag\\
\geq\; & \Big(2-\mu- 2 C_{h}\Big)\int\limits_S^T E_{\gamma,\delta}(t) dt
 -\frac{1}{2}\Big(2-\mu-2 C_{h}\Big)\int\limits_S^T\left(\gamma w^2(\ell)+\delta \psi^2(\ell)\right)dt
 \notag\\
&\hskip 7cm
 +2(F(T)-F(S))+\frac{\mu}{2}\int\limits_0^\ell \left( \rho w_t w +I_\rho \psi_t  \psi \right) dx|_S^T,
 \end{align*}
where we use \eqref{eq:hp_a} and \eqref{final_term}.
As we have seen in
Proposition \ref{F_estimate} and \ref{lower_order_estimate},
 the last three terms can be estimated by $E_{\gamma,\delta}(T)-E_{\gamma,\delta}(S)$ and, consequently, we 
obtain the estimate
\begin{align}\label{energy_estimate}
&\Big(2-\mu- 2 C_{h}\Big)\int\limits_S^T E_{\gamma,\delta}(t) dt \leq  \int\limits_S^T\ell  \left\{ \Big(\rho+\frac{\alpha^2}{K(\ell)}\Big)w_t^2(t,\ell)+\Big(I_\rho+\frac{\beta^2}{EI(\ell)}\Big)\psi_t^2(t,\ell) \right\} dt\notag\\
&  +\int\limits_S^T\left(\left(\frac{\ell}{K(\ell)}\gamma-\frac{\mu}{2}+\Big(2-\mu- 2C_{h}\Big)\right) \gamma w^2(t,\ell)+\left(\frac{\ell}{EI(\ell)}\delta-\frac{\mu}{2}+\Big(2-\mu-2 C_{h}\Big)\right)\delta\psi^2(t,\ell)\right) dt \notag\\
& \quad\quad +\int\limits_S^T\left(\alpha\left(\frac{2\ell\gamma}{K(\ell)}-\frac{\mu}{2}\right) w_t(t,\ell)w(t,\ell)+\delta\left(\frac{2\ell\delta}{EI(\ell)}-\frac{\mu}{2}\right)\psi_t(t,\ell) \psi(t,\ell)\right) dt\\
&\quad\quad\quad\quad  +\left(2 C_{F}+{\mu\over2} C_{NL} \right)(E_{\gamma,\delta}(T)+E_{\gamma,\delta}(S)).
\notag\end{align}
%
By applying the Cauchy-Schwartz inequality to the mixed terms and using the fact that the system is dissipative, we obtain
\begin{align}\label{energy_inequality}
&(2-\mu- C_{h})\int\limits_S^T E_{\gamma,\delta}(t) dt \leq  \int\limits_S^T  \left\{ \alpha\eta_{11}w_t^2(t,\ell)+\beta\eta_{12}\psi_t^2(t,\ell) \right\} dt
\\
& \hskip 4cm  +\int\limits_S^T\left(\eta_{21}\gamma w^2(t,\ell)+\eta_{22}\delta\psi^2(t,\ell)\right) dt +2\left(2 C_{F}+{\mu\over2} C_{NL} \right)E_{\gamma,\delta}(S), \notag
\end{align}
where   constants $\eta_{ij},\; i,j =1,2,$ are  positive.
 We recall that
$$\frac{d}{dt}E_{\gamma,\delta}(t)=-\alpha w_t(t,\ell)^2-\beta \psi_t(t,\ell)^2.$$
 Then \eqref{energy_inequality} turns into
\begin{align}\label{energy_inequality_1}
&(2-\mu-  C_{h})\int\limits_S^T E_{\gamma,\delta}(t) dt\leq
\int\limits_S^T\left(\eta_{21}\gamma\w(t,\ell)^2+\eta_{22}\delta\psi(t,\ell)^2\right) dt\notag\\
&\hskip 4cm  +\left(\max\{\eta_{11},\eta_{12}\}+2\left(2 C_{F}+{\mu\over2} C_{NL} \right)\right)E_{\gamma,\delta}(S).
\end{align}

\begin{remark}
In \eqref{energy_inequality_1}, we denote the constant in front of $E_{\gamma,\delta}(S)$ by $C'$. The goal is to estimate the integral $\int\limits_S^T\left(\eta_{212}\gamma\w(t,\ell)^2+\eta_{22}\delta\psi(t,\ell)^2\right) dt$ with respect to $E_{\gamma,\delta}(S)$. Clearly, we have that 
$$
\int\limits_S^T\left(\eta_{21}\gamma\w(t,\ell)^2+\eta_{22}\delta\psi(t,\ell)^2\right) dt \le  2\max\{\eta_{21},\eta_{22}\}\int\limits_S^ T E_{\gamma,\delta}(t)ds.
$$
We can let $\varepsilon := 2\max\{\eta_{12},\eta_{21}\}$ small enough such that
$2-\mu-  C_{h}-\varepsilon>0$, and consequently,  $(2-\mu-  C_{h}-\varepsilon)\int\limits_S^T E_{\gamma,\delta}(t) dt\leq C'E_{\gamma,\delta}(S)$.
 However, that is too coarse. Indeed, similar to \cite{AlabauCannarsaLeugering2017}, we proceed differently to prove the exponential stability.
\end{remark}

Consider the following auxiliary degenerated elliptic problem
\begin{align}\label{elliptic}
& (K(x)(z_x+\xi))_x=0,\; x\in (0,\ell)\notag\\
& (EI(x)\xi_x)_x -K(x)(z_x+\xi)=0,\; x\in (0,\ell),\notag\\
& z(0)=0, \; \text{ if } \mu_K\in[0,1),\; \lim\limits_{x\rightarrow 0^+}K(x)(z_x+\xi)(x)=0\;  \text{ if } \; \mu_K\in [1,2),\\
& \xi(0)=0, \; \text{ if } \mu_{EI}\in[0,1),\; \lim\limits_{x\rightarrow 0^+}EI(x)(z_x+\xi)(x)=0\;  \text{ if } \; \mu_EI\in [1,2),\notag\\
& K(\ell)(z_x(\ell)+\xi(\ell))+\gamma z(\ell)=\lambda, \; EI(\ell)\xi(\ell)+\delta \xi(\ell)=\sigma.\notag
\end{align}
We notice that problem \eqref{elliptic} can be formulated in the variational form. To this end, we recall the bilinear form \eqref{bilinear_gd}
and introduce the linear functional $
L(u,\chi):=\lambda u(\ell)+\sigma \chi(\ell)$. We pose the problem
\begin{equation}\label{elliptic_var}
B_{\gamma,\delta}((z,\xi),(u,\chi))= L(u,\chi), \;\; \forall (u,\chi) \in ({\cal V}_{R,w/s})^2,\end{equation}
It is now clear, according to the Lax-Milgram theorem, that \eqref{elliptic_var} admits a unique solution $(z,\xi)$. Note that   $
\|| (z,\xi) \||_{\gamma,\delta}^2= B_{\gamma,\delta}((z,\xi),(z,\xi)) \geq \gamma z^2(\ell)+\delta \xi^2(\ell)$ holds. Furthermore, by \eqref{elliptic_var},
\begin{equation}
\label{lambda}
\|| (z,\xi) \||_{\gamma,\delta}^2 \leq
{|\lambda| \over\sqrt{\gamma}} \sqrt{\gamma}|u(\ell)| + {|\sigma| \over\sqrt{\delta}} \sqrt{\delta}|\chi(\ell)| \le \max\Big\{\frac{1}{\sqrt{\gamma}},\frac{1}{\sqrt{\delta}}\Big\}\left( |\lambda|+|\sigma|\right)\|| (z,\xi) \||_{\gamma,\delta}.
\end{equation}
Therefore, the solution to \eqref{elliptic} satisfies
\begin{equation}\label{lambdasigma}
\gamma z(\ell)^2+\delta \xi(\ell)^2\leq \|| (z,\xi)\||_{\gamma,\delta}^2 \le 
2\max\Big\{\frac{1}{\gamma},\frac{1}{\delta}\Big\}\left( \lambda^2+\sigma^2\right).
\end{equation}

We take  $\lambda=w(t,\ell), \sigma=\psi(t,\ell)$, such that $(z,\xi)$ is now quasi-static, as it depends on $t$ through $\lambda, \sigma$ and multiply the system \eqref{Timo_weak_strong_Neumann} with \eqref{Timo_feedback} by $z$ and $\xi$  and integrate by parts.
\begin{align}\label{zxi}
  0=& \int\limits_S^T \int\limits_0^\ell \left\{ \rho w_{tt} z +I_\rho \psi_{tt}\xi-(K(w_x+\psi))_xz-(EI\psi_x)_x\xi+K(w_x+\psi)\xi\right\}dxdxt\notag\\
 =& \int\limits_0^\ell( \rho w_t z +I_\rho \psi_{t} \xi )dx |_S^T-\int\limits_S^T\int\limits_0^\ell \big( \rho w_tz_t+I_\rho \psi_t\xi_t \big)dxdt \\ & \qquad\qquad +\int\limits_S^T (\alpha w_t(t,\ell)z(t,\ell)+\beta \psi_t(t,\ell) \xi(t,\ell)) dt +  \int\limits_S^T  (w^2(t,\ell)+\psi^2(t,\ell) ) dt.\notag\end{align}
 From \eqref{lambdasigma}, we have 
 \begin{equation}
 \label{zxibc}
z^2(t,\ell)+\xi^2(t,\ell)\leq 2\max\Big\{\frac{1}{\gamma^2},\frac{1}{\delta^2}\Big\}(w^2(t,\ell)+\psi^2(t,\ell)).
 \end{equation}
We may now take the time derivative of \eqref{elliptic} and obtain the same estimation as \eqref{lambdasigma} for $(z_t,\xi_t)$ corresponds to the data $w_t(t,\ell), \psi_t(\ell)$.  Combining these with Proposition \ref{Poincare_gd}, we obtain
\begin{equation}\label{zxi_L2}
    \int\limits_0^\ell\big( z_t^2(t,x)+\xi_t^2(t,x)\big)dx \leq 2 \max\Big\{\frac{1}{\gamma^2},\frac{1}{\delta^2}\Big\}
    C_{N,K,EI}\big(w_t^2(t,\ell)+\psi_t^2(t,\ell)\big).
\end{equation}
Therefore, using the Cauchy-Schwartz inequality and \eqref{zxi_L2}, we obtain for some $\omega>0$\begin{align}
   & \int\limits_S^T\int\limits_0^\ell \big(\rho  w_t z_t+I_\rho \psi_t\xi_t\big) dx dt\leq  \frac{\omega}{2}\int\limits_S^T\int\limits_0^\ell \big(\rho  w_t^2 +I_\rho \psi_t^2\big) dx dt+\frac{1}{2\omega}\max\{\rho, I_\rho\}\int\limits_S^T\int\limits_0^\ell \big( z_t^2 + \xi_t^2\big) dx dt\notag\\
   \leq\; & \omega \int\limits_S^T E_{\gamma,\delta}(t)ds+ \frac{1}{\omega}\max\{\rho, I_\rho\}\max\Big\{\frac{1}{\gamma^2},\frac{1}{\delta^2}\Big\}
   \max\Big\{\frac{1}{\alpha},\frac{1}{\beta}\Big\}C_{N,K,EI}
   (E_{\gamma,\delta}(S)-E_{\gamma,\delta}(T)) \notag\\
  =: \; & \omega \int\limits_S^T E_{\gamma,\delta}(t)ds+ \frac{1}{\omega}C''(E_{\gamma,\delta}(S)-E_{\gamma,\delta}(T))\leq
   \omega \int\limits_S^T E_{\gamma,\delta}(t)ds+ \frac{1}{\omega}C''E_{\gamma,\delta}(S), \end{align}
   where
   $$
   C'' := \max\{\rho, I_\rho\}\max\Big\{\frac{1}{\gamma^2},\frac{1}{\delta^2}\Big\}
   \max\Big\{\frac{1}{\alpha},\frac{1}{\beta}\Big\}C_{N,K,EI}.
   $$
Now, using again the Cauchy-Schwartz inequality and \eqref{lambdasigma}, we have
\begin{align}
&\int\limits_0^\ell (\rho w_t z +I_\rho \psi_{t} \xi )dx\leq \frac{1}{2}\int\limits_0^\ell( \rho w_t^2+I_\rho \psi_t^2 )dx +
\frac{1}{2}\int\limits_0^\ell( \rho z^2+I_\rho \xi^2 )dx \notag\\
\leq\;& \frac{1}{2}\int\limits_0^\ell( \rho w_t^2+I_\rho \psi_t^2) dx +
 \max\{\rho,I_\rho\}\max\Big\{\frac{1}{\gamma^3},\frac{1}{\delta^3}\Big\} \big(\gamma w(t,\ell)^2+\delta\psi(t,\ell)^2 \big)\notag\\
  \leq \;  &\Big(1+ 2\max\{\rho,I_\rho\}\max\Big\{\frac{1}{\gamma^3},\frac{1}{\delta^3}\Big\} \Big)E_{\gamma,\delta}(t) 
=: C''' E_{\gamma,\delta}(t).
\end{align}
where
$$
C''' := 1+ 2\max\{\rho,I_\rho\}\max\Big\{\frac{1}{\gamma^3},\frac{1}{\delta^3}\Big\}.
$$
%
Moreover, from \eqref{dissipative}, \eqref{lambdasigma}, we have
 \begin{align}
& \int\limits_S^T \big(\alpha w_t(t,\ell)z(t,\ell)+ \beta\psi_t(t,\ell)\xi(t,\ell)\big) dt
\notag \\
\leq\;& \int\limits_S^T \Big(\frac{\alpha}{2\omega} w_t(t,\ell)^2+ \frac{\alpha\omega}{2}z(t,\ell)^2+\frac{\beta}{2\omega} \psi_t(t,\ell)^2 +\frac{\omega\beta}{2}\xi(t,\ell)^2 \Big)dt\notag\\
\leq \;& \frac{1}{2\omega}\int\limits_S^T \big(\alpha w_t(t,\ell)^2+\beta \psi_t(t,\ell)^2\big) dt+ \frac{\omega}{2}\max\Big\{\frac{\alpha}{\gamma},\frac{\beta}{\delta}\Big\} \int\limits_S^T(\gamma z(t,\ell)^2+\delta\xi(t,\ell)^2)dt\notag\\
\leq \; & \frac{1}{\omega}E(S) +\omega\max\Big\{\frac{\alpha}{\gamma^3},\frac{\beta}{\delta^3}\Big\} \int\limits_S^T E_{\gamma,\delta}(t)dt .
 \end{align}
 Putting things together, we finally obtain from \eqref{zxi} that
 \begin{align}\label{wpsiatell}
  &\int\limits_S^T \big(\eta_{21} w(t,\ell)^2+\eta_{22}\psi(t,\ell)^2 \big)dt \leq \max\big\{\eta_{21},\eta_{22}\big\} \Big[\omega \Big( 1+\max\Big\{\frac{\alpha}{\gamma^3},\frac{\beta}{\delta^3}\Big\} \Big) \int\limits_S^T E_{\gamma,\delta}(t)dt\notag\\   &\qquad\qquad \hskip 6cm +\frac{1}{\omega}(1+ C'' )E_{\gamma,\delta}(S)
 +  2C'''E_{\gamma,\delta}(S) \Big].
 \end{align}
 Hence, \eqref{energy_inequality_1} becomes
\begin{align}\label{tildec}
&\Big(2-\mu-  C_{h} - \omega\max\big\{\eta_{21},\eta_{22}\big\}\Big( 1+2\max\Big\{\frac{\alpha}{\gamma^3},\frac{\beta}{\delta^3}\Big\}  \Big) \Big)\int\limits_S^T E_{\gamma,\delta}(t)dt \leq \tilde{C}(\omega) E_{\gamma,\delta}(S).
\end{align}
where
$$
\tilde{C}(\omega):= \frac{1}{\omega}\max\big\{\eta_{21},\eta_{22}\big\}(1+ C'')+2\max\big\{\eta_{21},\eta_{22}\big\}C''' +\max\{\eta_{11},\eta_{12}\}+2\left(2 C_{F}+{\mu\over2} C_{NL} \right).
$$
We now choose $\omega$ as
\begin{align}\label{omega}
\omega:= \frac{2-\mu-C_{h}}{2\max\big\{\eta_{21},\eta_{22}\big\}\Big( 1+2\max\Big\{\frac{\alpha}{\gamma^3},\frac{\beta}{\delta^3}\Big\}  \Big)}
\end{align}
 and insert \eqref{omega} into $\eqref{tildec}$. This results in
 \begin{equation}\label{energy_inequality_final}
 \Big(2-\mu-  C_{h}\Big)\int\limits_S^T E_{\gamma,\delta}(t)dt\leq 2 \tilde{C}(\omega) E(S).
 \end{equation}
We define $$\kappa:=\frac{2-\mu- C_{h}}{2\tilde{C}(\omega)}.$$ By Theorem 8.1 of \cite{Komornik1989}, \eqref{energy_inequality_final} implies
\begin{theorem}
Let \eqref{small} be satisfied. Then the solutions of \eqref{Timo_feedback} satisfy
    \begin{equation}
     E_{\gamma,\delta}(t)\leq E_{\gamma,\delta}(0)e^{1-\kappa t}, \;\; t\geq \frac{1}{\kappa}.
 \end{equation}
\end{theorem}

\begin{remark}
As stated in Remark \label{smallh}, we can choose $h$ sufficiently small such that $C_h$ becomes small. Therefore, condition \eqref{small} is reasonable. On the other hand, when $\mu \to2$, i.e., $\max\{\mu_K,\; \mu_{EI}\| \to2, $
he questions of observability, exact controllability, and exponential stability remain open.
 \end{remark}

\end{document}